\documentclass[twoside,draft,12pt]{article}
\usepackage{amssymb}

\usepackage{makeidx}
\usepackage{a4wide}
\usepackage{amsthm}
\usepackage{amsmath}
\usepackage[all]{xy}
\usepackage{enumerate}
\usepackage{fancyhdr}

\newtheorem{theorem}{Theorem}[section]
\newtheorem{proposition}[theorem]{Proposition}
\newtheorem{corollary}[theorem]{Corollary}
\newtheorem{lemma}[theorem]{Lemma}
\theoremstyle{definition}
\newtheorem{example}[theorem]{Example}
\newtheorem{c-example}[theorem]{Counter Example}
\newtheorem*{Notation}{Notation}
\newtheorem*{Beweis}{Proof}
\newtheorem{definition}[theorem]{Definition}
\newtheorem{punto}[theorem]{}
\theoremstyle{remark}
\newtheorem{remark}[theorem]{Remark}
\newtheorem{remarks}[theorem]{Remarks}
\CompileMatrices
\input{tcilatex}

\begin{document}

\title{On Coreflexive Coalgebras and Comodules over Commutative Rings\thanks{%
MSC (2000): 16D90, 16W30, 16Exx \newline
Keywords: Dual Coalgebras, Dual Comodules, Coreflexive Coalgebras, Reflexive
Algebras, Coreflexive Comodules, Reflexive Modules, Linear Weak Topology}}
\author{\textbf{Jawad Y. Abuhlail} \\
%EndAName
Mathematics Department\\
Birzeit University\\
P.O.Box 14, Birzeit - Palestine \\
jabuhlail@birzeit.edu}
\date{}
\maketitle

\begin{abstract}
In this note we study dual coalgebras of algebras over arbitrary
(noetherian) commutative rings. We present and study a generalized notion of
coreflexive comodules and use the results obtained for them to characterize
the so called coreflexive coalgebras. Our approach in this note is an
algebraically topological one.
\end{abstract}

\subsection*{Introduction}

The concept of \emph{coreflexive coalgebras} was studied, in the case of
commutative base fields, by several authors. An algebraic approach was
presented by E. Taft (\cite{Taf72}, \cite{Taf77}), while a topological one
was presented mainly by D. Radford (\cite{HR74}, \cite{Rad73}) and studied
by several authors (e.g. \cite{Miy75}, \cite{Wit79}). In this note we
present and study a generalized concept of \emph{coreflexive comodules }and
use it to characterize coreflexive coalgebras over commutative (noetherian)
rings. In particular we generalize results from the papers mentioned above
from the case of base fields to the case of arbitrary (noetherian)
commutative ground rings.

Throughout this paper $R$ denotes a commutative ring with $1_{R}\neq 0_{R}.$
We consider $R$ as a left and a right linear topological ring with the \emph{%
discrete topology. }The category of $R$-(bi)modules will be denoted by $%
\mathcal{M}_{R}.$ The unadorned $-\otimes -$ and $\mathrm{Hom}$ mean $%
-\otimes _{R}-$ and $\mathrm{Hom}_{R}$ respectively. For an $R$-module $M,$
an $R$-submodule $K\subset M$ will be called $N$\emph{-pure} for some $R$%
-module $N,$ if the canonical $R$-linear mapping $\iota _{K}\otimes
id_{N}:K\otimes _{R}N\rightarrow M\otimes _{R}N$ is injective. We call $%
K\subset M$ \emph{pure }(in the sense of Cohn), if it's $N$-pure for every $%
R $-module $N.$ For every $R$-module $L,$ we denote with $L^{\ast }$ the 
\emph{algebraic dual }$R$-module of all $R$-linear maps from $L$ to $R.$

Let $S$ be a ring. We consider every left (right) $S$-module $K$ as a right
(a left) module over $\mathrm{End}(_{S}K)^{op}$ ($\mathrm{End}(K_{S})$) and
as a left (a right) module over $\mathrm{Biend}(_{S}K):=\mathrm{End}(K_{%
\mathrm{End}(_{S}K)^{op}})$ ($\mathrm{Biend}(K_{S}):=\mathrm{End}(_{\mathrm{%
End}(K_{S})}K)^{op}$), the \emph{ring of biendomorphisms} of $K$ (e.g. 
\cite[6.4]{Wis88}).

Let $A$ be an $R$-algebra and $M$ be an $A$-module. An $A$-submodule $%
N\subset M$ will be called $R$\emph{-cofinite}, if $M/N$ is finitely
generated in $\mathcal{M}_{R}.$ The class of all $R$-cofinite $A$-submodules
of $M$ is denoted with $\mathcal{K}_{M}.$ We call $M$ \emph{cofinitely }$R$%
\emph{-cogenerated}, if $M/N$ is $R$-cogenerated for every $R$-cofinite $A$%
-submodule. With $\mathcal{K}_{A}$ we denote the class of \emph{all }$R$%
-cofinite $A$-ideals and define 
\begin{equation*}
A^{\circ }:=\{f\in A^{\ast }\mid f(I)=0\text{ for some }R\text{-cofinite
ideal }I\vartriangleleft A\}.
\end{equation*}
If $\mathcal{K}_{A}$ is a filter (e.g. $R$ is a noetherian ring), then $%
A^{\circ }\subseteq A^{\ast }$ is an $R$-submodule with equality, iff $_{R}A$
is finitely generated projective.

We assume the reader is familiar with the theory of Hopf Algebras. For any
needed definitions or results the reader may refer to any of the classical
books on the subject (e.g. \cite{Swe69}, \cite{Abe80} and \cite{Mon93}). For
an $R$-coalgebra $(C,\Delta _{C},\varepsilon _{C})$ and an $R$-algebra $%
(A,\mu _{A},\eta _{A})$ we consider $\mathrm{Hom}_{R}(C,A)\;$as an $R$%
-algebra with multiplication the \emph{convolution product }$(f\star
g)(c):=\sum f(c_{1})g(c_{2})$ and unity $\eta _{A}\circ \varepsilon _{C}.$

\section{Preliminaries}

In this section we present some definitions and lemmata.

\begin{definition}
Let $(C,\Delta _{C},\varepsilon _{C})$ be an $R$-coalgebra. We call an $R$%
-submodule $K\subset C:$

an $R$\emph{-subcoalgebra}, iff $K\subset C$ is pure and $\Delta
_{C}(K)\subset K\otimes _{R}K;$

a $C$\emph{-coideal}, iff $K\subset \mathrm{Ke}(\varepsilon _{C})$ and 
\begin{equation*}
\Delta _{C}(K)\subset \mathrm{Im}(\iota _{K}\otimes id_{C})+\mathrm{Im}%
(id_{C}\otimes \iota _{K});
\end{equation*}

a \emph{right }$C$\emph{-coideal} (resp. a \emph{left }$C$\emph{-coideal}, a 
$C$\emph{-bicoideal}), if $K\subset C$ is $C$-pure and $\Delta
_{C}(K)\subset K\otimes _{R}C$ (resp. $\Delta _{C}(K)\subset C\otimes _{R}K,$
$\Delta _{C}(K)\subset (K\otimes _{R}C)\cap (C\otimes _{R}K)$).
\end{definition}

\begin{punto}
\textbf{Subgenerators}. Let $A$ be an $R$-algebra and $K$ be a left $A$%
-module. We say a left $A$-module $N$ is $K$\emph{-subgenerated}, if $N$ is
isomorphic to a submodule of a $K$-generated left $A$-module (equivalently,
if $N$ is kernel of a morphism between $K$-generated left $A$-modules). The 
\emph{full }subcategory of $_{A}\mathcal{M},$ whose objects are the $K$%
-subgenerated left $A$-modules is denoted by $\sigma \lbrack _{A}K].$ In
fact $\sigma \lbrack _{A}K]\subseteq $ $_{A}\mathcal{M}$ is the \emph{%
smallest }Grothendieck full subcategory that contains $K.$ If $M$ is a left $%
A$-module, then 
\begin{equation*}
\mathrm{Sp}(\sigma \lbrack _{A}K],M):=\sum \{f(N)\mid \text{ }f\in \mathrm{%
Hom}_{A-}(N,M),\text{ }N\in \sigma \lbrack _{A}K]\}
\end{equation*}
is the biggest $A$-submodule of $M$ that belongs to $\sigma \lbrack _{A}K].$
The reader is referred to \cite{Wis88} and \cite{Wis96} for the well
developed theory of categories of this type.
\end{punto}

\subsection*{The linear weak topology}

\begin{punto}
$R$\textbf{-Pairings.} An $R$\emph{-pairing} $P=(V,W)$ consists of $R$%
-modules $V,W$ with an $R$-bilinear form 
\begin{equation*}
\alpha :V\times W\rightarrow R,\text{ }(v,w)\mapsto <v,w>.
\end{equation*}
If the induced $R$-linear mapping $\kappa _{P}:V\rightarrow W^{\ast }$
(resp. $\chi _{P}:W\rightarrow V^{\ast }$) is injective, then we call $P$ 
\emph{left non-degenerating}\textbf{\ (}resp. \emph{right non-degenerating}%
). If both $\kappa _{P}$ and $\chi _{P}$ are injective, then we call $P$ 
\emph{non-degenerating.}

For $R$-pairings $(V,W)$ and $(V^{\prime },W^{\prime })$ a morphism 
\begin{equation*}
(\xi ,\theta ):(V^{\prime },W^{\prime })\rightarrow (V,W)
\end{equation*}
consists of $R$-linear mappings $\xi :V\rightarrow V^{\prime }$ and $\theta
:W^{\prime }\rightarrow W,$ such that 
\begin{equation*}
<\xi (v),w^{\prime }>=<v,\theta (w^{\prime })>\text{ for all }v\in V\text{
and }w^{\prime }\in W^{\prime }.
\end{equation*}
The $R$-pairings with the morphisms described above (and the usual
composition of pairings) build a category which we denote with $\mathcal{P}.$
If $P=(V,W)$ is an $R$-pairing, $V^{\prime }\subset V$ is an $R$-submodule
and $W^{\prime }\subset W$ is a (pure) $R$-submodule with $<V^{\prime
},W^{\prime }>=0,$ then $Q:=(V/V^{\prime },W^{\prime })$ is an $R$-pairing, $%
(\pi ,\iota _{K}):(V/V^{\prime },W^{\prime })\rightarrow (V,W)$ is a
morphism in $\mathcal{P}$ and we call $Q\subset P$ a (\emph{pure}) $R$\emph{%
-subpairing}.
\end{punto}

\begin{Notation}
Let $P=(V,W)$ be an $R$-pairing. For subsets $X\subseteq V$ resp. $%
K\subseteq W$ set 
\begin{equation*}
X^{\bot }:=\{w\in W|\text{ }<X,w>=0\}\text{ resp. }K^{\bot }:=\{v\in V|\text{
}<v,K>=0\}.
\end{equation*}
We say $X\subseteq V$ (resp. $K\subseteq W$) is \emph{orthogonally closed}
w.r.t. $P,$ if $X=X^{\bot \bot }$ (resp. $K=K^{\bot }{}^{\bot }$). In case $%
V=W^{\ast },$ then we set for every subset $X\subset W^{\ast }$ (resp. $%
K\subset W$) $\mathrm{Ke}(X)=X^{\bot }$ (resp. $\mathrm{An}(K)=:K^{\bot }$).
\end{Notation}

\begin{punto}
\label{lst-CA}{\normalsize \ }Let $P=(V,W)$ be an $R$-pairing. Then the
class of $R$-submodules of $V:$%
\begin{equation*}
\mathcal{F}{\normalsize (0_{V}):=\{K^{\bot }|}\text{ }K\subset W\text{ is a
finitely generated }R\text{-submodule}\}
\end{equation*}
is a filter basis consisting of $R$-submodule of $V$ and induces on $V$ a
topology, the so called \emph{linear weak topology} $V[{\frak{T}}_{ls}(W)],$
such that $(V,V[\frak{T}_{ls}(W)])$ is a linear topological right $R$-module
and $\mathcal{F}(0_{V})$ is a neighbourhood basis of $0_{V}.$ In particular
we call $W^{\ast }[{\frak{T}}_{ls}(W)]$ the finite topology. The properties
of this topology were studied by several authors in the case of commutative
base fields (e.g. \cite{Kot66}, \cite{KN63}, \cite{Rad73}). We refer mainly
to the recent work of the author \cite{Abu} for the case of arbitrary ground
rings.
\end{punto}

\qquad

\subsection*{The $\protect\alpha$-condition}

In a joint work with J. G\'{o}mez-Torrecillas and J. Lobillo \cite{AG-TL2001}
on the category of comodules of coalgebras over arbitrary commutative base
rings, we presented the so called $\alpha $\emph{-condition.} That condition
has shown to be a natural assumption in the author's study of \emph{duality
theorems }for Hopf algebras \cite{Abu2001}. We refer mainly to \cite{Abu}
for the properties of the such pairings over arbitrary ground rings.

\begin{punto}
$R$\textbf{-Pairings.}\label{P-alph} We say an $R$-pairing $P=(V,W)$
satisfies the $\alpha $\emph{-condition}\textbf{\ }(or $P$ is an $\alpha $%
\emph{-pairing}), if for every $R$-module $M$ the following map is injective 
\begin{equation}
\alpha _{M}^{P}:M\otimes _{R}W\rightarrow \mathrm{Hom}_{R}(V,M),\text{ }\sum
m_{i}\otimes w_{i}\mapsto \lbrack v\mapsto \sum m_{i}<v,w_{i}>].  \label{alp}
\end{equation}

With $\mathcal{P}^{\alpha }\subset \mathcal{P}$ we denote the \emph{full}
subcategory of $R$-pairings satisfying the $\alpha $-condition. We call an $%
R $-pairing $P=(V,W)$ \emph{dense}, if $\kappa _{P}(V)\subseteq W^{\ast }$
is dense (considering $W^{\ast }$ with the finite topology). It's easy to
see that $\mathcal{P}^{\alpha }\subset \mathcal{P}$ is closed under pure $R$%
-subpairings.

We say an $R$-module $W$ \emph{satisfies the }$\alpha $\emph{-condition}, if
the $R$-pairing $(W^{\ast },W)$ satisfies the $\alpha $-condition, i.e. for
every $R$-module $M$ the canonical $R$-linear morphism $\alpha
_{M}^{W}:M\otimes _{R}W\rightarrow \mathrm{Hom}_{R}(W^{\ast },M)$ in
injective (equivalently, if $_{R}W$ is \emph{locally projective} in the
sense of B. Zimmermann-Huisgen \cite{ZH76}).
\end{punto}

\begin{remark}
\label{flat}\cite[Remark 2.2]{Abu} Let $P=(V,W)\in \mathcal{P}^{\alpha }.$
Then $_{R}W$ is $R$-cogenerated and flat. If $R$ is perfect, then $_{R}W$
turns to be projective.
\end{remark}

\qquad\ 

\begin{Notation}
Let $W,W^{\prime }$ be $R$-modules and consider for any $R$-submodules $%
X\subseteq W^{\ast }$ and $X^{\prime }\subseteq W^{\prime \ast }$ the
canonical $R$-linear mapping 
\begin{equation*}
\delta :X\otimes _{R}X^{\prime }\rightarrow (W\otimes _{R}W^{\prime })^{\ast
}.
\end{equation*}
For $f\in X$ and $g\in X^{\prime }$ set $f\underline{\otimes }g=\delta
(f\otimes g),$ i.e. 
\begin{equation*}
(f\underline{\otimes }g)(\sum w_{i}\otimes w_{i}^{\prime }):=\sum
f(w_{i})g(w_{i}^{\prime })\text{ for every }\sum w_{i}\otimes w_{i}^{\prime
}\in W\otimes _{R}W^{\prime }.
\end{equation*}
\end{Notation}

\section{Measuring $R$-pairings}

\begin{punto}
\label{C-bimod}For an $R$-coalgebra $C$ and an $R$-algebra $A$ we call an $R$%
-pairing $P=(A,C)$ a \emph{measuring }$R$\emph{-pairing}, if the induced
mapping $\kappa _{P}:A\rightarrow C^{\ast }$ is an $R$-algebra morphism. In
this case $C$ is an $A$-bimodule through the left and the right $A$-actions 
\begin{equation}
a\rightharpoonup c:=\sum c_{1}<a,c_{2}>\text{ and }c\leftharpoonup a:=\sum
<a,c_{1}>c_{2}\text{ for all }a\in A,\text{ }c\in C.  \label{C-r}
\end{equation}
Let $(A,C)$ and $(B,D)$ be measuring $R$-pairings. We say a morphism of $R$%
-pairings $(\xi ,\theta ):(B,D)\rightarrow (A,C)$ is a \emph{morphism of
measuring }$R$\emph{-pairings,} if $\xi :A\rightarrow B$ is an $R$-algebra
morphism and $\theta :D\rightarrow C$ is an $R$-coalgebra morphism. The
category of measuring $R$-pairings and morphisms described above will be
denoted by $\mathcal{P}_{m}.$ With $\mathcal{P}_{m}^{\alpha }\subset 
\mathcal{P}_{m}$ we denote the \emph{full }subcategory of measuring $R$%
-pairings satisfying the $\alpha $-condition (we call these \emph{measuring }%
$\alpha $\emph{-pairings}). If $P=(A,C)$ is a measuring $R$-pairing, $%
D\subset C$ is an $R$-subcoalgebra and $I\vartriangleleft A$ is an ideal
with $<I,D>=0,$ then $Q:=(A/I,C)$ is a measuring $R$-pairing, $(\pi
_{I},\iota _{D}):(A/I,D)\rightarrow (A,C)$ is a morphism in $\mathcal{P}_{m}$
and we call $Q\subset P$ a \emph{measuring }$R$\emph{-subpairing}. Since by
convention an $R$-coalgebra is a pure $R$-submodule, it's easy to see that $%
\mathcal{P}_{m}^{\alpha }\subset \mathcal{P}_{m}$ is closed against
measuring $R$-subpairings.
\end{punto}

\begin{lemma}
\label{prop-mes}Let $P=(A,C),$ $Q=(B,D)\in \mathcal{P}_{m}$ and $(\xi
,\theta ):(B,D)\rightarrow (A,C)$ be a morphism of $R$-pairings.

\begin{enumerate}
\item  Assume that $P\otimes P:=(A\otimes _{R}A,C\otimes _{R}C)$ is right
non-degenerating \emph{(}i.e. $\chi :=\chi _{P\otimes P}:C\otimes
_{R}C\hookrightarrow (A\otimes _{R}A)^{\ast }$ is an embedding\emph{)}. If $%
\xi $ is an $R$-algebra morphism, then $\theta $ is an $R$-coalgebra
morphism. If $A$ is commutative, then $C$ is cocommutative.

\item  If $Q$ is left non-degenerating \emph{(}i.e. $B\overset{\kappa _{Q}}{%
\hookrightarrow }D^{\ast }$ is an embedding\emph{)} and $\theta $ is an $R$%
-coalgebra morphism, then $\xi $ is an $R$-algebra morphism. If $C$ is
cocommutative and $P$ is left non-degenerating \emph{(}i.e.\emph{\ }$%
A\subseteq C^{\ast }$\emph{)}, then $A$ is commutative.
\end{enumerate}
\end{lemma}

\begin{Beweis}
\begin{enumerate}
\item  If $\xi $ is an $R$-algebra morphism, then we have for arbitrary $%
d\in D,$ $a,\widetilde{a}\in A:$%
\begin{equation*}
\begin{tabular}{lll}
$\chi (\sum \theta (d)_{1}\otimes \theta (d)_{2})(a\otimes \widetilde{a})$ & 
$=$ & $\sum <a,\theta (d)_{1}><\widetilde{a},\theta (d)_{2}>$ \\ 
& $=$ & $<a\widetilde{a},\theta (d)>$ \\ 
& $=$ & $<\xi (a\widetilde{a}),d>$ \\ 
& $=$ & $<\xi (a)\xi (\widetilde{a}),d>$ \\ 
& $=$ & $\sum <\xi (a),d_{1}><\xi (\widetilde{a}),d_{2}>$ \\ 
& $=$ & $\sum <a,\theta (d_{1})><\widetilde{a},\theta (d_{2})>$ \\ 
& $=$ & $\chi (\sum \theta (d_{1})\otimes \theta (d_{2}))(a\otimes 
\widetilde{a}).$%
\end{tabular}
\end{equation*}
By assumption $\chi $ is injective and so $\sum \theta (d)_{1}\otimes \theta
(d)_{2}=\sum \theta (d_{1})\otimes \theta (d_{2})$ for every $d\in D,$ i.e. $%
\theta $ is an $R$-coalgebra morphism.

If $A$ is commutative, then we have for all $c\in C$ and $a,\widetilde{a}\in
A:$%
\begin{equation*}
\begin{tabular}{lllll}
$\chi (\sum c_{1}\otimes c_{2})(a\otimes \widetilde{a})$ & $=$ & $\sum
<a,c_{1}><\widetilde{a},c_{2}>$ & $=$ & $<a\widetilde{a},c>$ \\ 
& $=$ & $<\widetilde{a}a,c>$ & $=$ & $\sum <a,c_{2}><\widetilde{a},c_{1}>$
\\ 
& $=$ & $\chi (\sum c_{2}\otimes c_{1})(a\otimes \widetilde{a}).$ &  & 
\end{tabular}
\end{equation*}
By assumption $\chi $ is injective and so $\sum c_{1}\otimes c_{2}=\sum
c_{2}\otimes c_{1}$ for every $c\in C,$ i.e. $C$ is cocommutative.

\item  Analogous to (1).$\blacksquare $
\end{enumerate}
\end{Beweis}

\begin{Notation}
Let $A$ be an $R$-algebra and $N$ be a left $A$-module (resp. a right $A$%
-module). For subsets $X,Y\subset N$ we set 
\begin{equation*}
(Y:X):=\{a\in A|\text{ }aX\subset Y\}\text{ (resp. }(Y:X):=\{a\in A|\text{ }%
Xa\subset Y\}\text{).}
\end{equation*}
If $Y=\{0_{N}\},$ then we set also $\mathrm{An}_{A}(X):=(0_{N}:X).\mathrm{\ }
$If $N$ is an $A$-bimodule, then we set for every subset $X\subset N:$%
\begin{equation*}
\mathrm{An}_{A}^{l}(X):=\{a\in A|\text{ }aX=0_{N}\}\text{ and }\mathrm{An}%
_{A}^{r}(X):=\{a\in A|\text{ }Xa=0_{N}\}.
\end{equation*}
\end{Notation}

\begin{punto}
\textbf{The }$C$-\textbf{adic topology.}\label{C-ad} Let $(A,C)\in \mathrm{P}%
_{m}$ and consider $C$ as a left $A$-module with the left $A$-action ``$%
\rightharpoonup $'' in (\ref{C-r}). Then the class of left $A$-ideals 
\begin{equation*}
\mathcal{B}{\normalsize _{C-}}(0_{A}):=\{\mathrm{An}_{A}^{l}(W)=(0:W)|\text{ 
}W=\{c_{1},...,c_{k}\}\subset C\text{ a finite subset}\}
\end{equation*}
is a neighbourhood basis of $0_{A}$ and induces on $A$ a topology, the so
called \emph{left }$C$\emph{-adic topology} $\mathcal{T}_{C-}(A),$ so that $%
(A,\mathcal{T}_{C-}(A))$ is a left linear topological $R$-algebra (see \cite
{AW97}, \cite{Ber94}). A left $A$-ideal $I\vartriangleleft _{l}A$ is open
w.r.t. $\mathcal{T}_{C-}(A),$ iff $A/I$ is $C$-subgenerated. If $\frak{T}$
is a left linear topology on $A,$ then the category of discrete left $(A,%
\frak{T})$-modules is equal to the category of $C$-subgenerated left $A$%
-modules $\sigma \lbrack _{A}C],$ iff $\frak{T}=${\normalsize $\mathcal{T}%
_{C-}(A).$} In particular we have for every left $A$-module $N:$%
\begin{equation*}
\mathrm{Sp}(\sigma \lbrack _{A}C],N)=\{n\in N|\text{ }\exists \text{ }%
F=\{c_{1},...,c_{k}\}\subset C\text{ with }(0_{C}:F)\subset (0_{N}:n)\}.
\end{equation*}
Analogously $C_{A}$ induces on $A$ a topology, the so called \emph{right }$C$%
\emph{-adic topology} $\mathcal{T}_{-C}(A),$ such that $(A,\mathcal{T}%
_{-C}(A))$ is a right linear topological $R$-algebra.
\end{punto}

\subsection*{Rational modules}

\begin{punto}
\label{rat-dar}Let $P=(A,C)$ be a measuring $\alpha $-pairing. Let $M$ be a
left $A$-module, $\rho _{M}:M\rightarrow \mathrm{Hom}_{R}(A,M)$ be the
canonical $A$-linear mapping and $\mathrm{Rat}^{C}(_{A}M):=\rho
_{M}^{-1}(M\otimes _{R}C).$ In case $\mathrm{Rat}^{C}(_{A}M)=M$ we call $M$
a $C$\emph{-rational }left $A$-module and define 
\begin{equation*}
\varrho _{M}:=(\alpha _{M}^{P})^{-1}\circ \rho _{M}:M\rightarrow M\otimes
_{R}C.
\end{equation*}
Analogously one defines the $C$\emph{-rational right }$A$\emph{-modules}.
With $\mathrm{Rat}^{C}(_{A}\mathcal{M})\subset $ $_{A}\mathcal{M}$ (resp. $%
^{C}\mathrm{Rat}(\mathcal{M}_{A})\subset \mathcal{M}_{A}$) we denote the 
\emph{full} subcategory of $C$-rational left (resp. right) $A$-modules.
\end{punto}

\begin{lemma}
\label{clos}\emph{(\cite[Lemma 2.2.7]{Abu2001})}\ Let $P=(A,C)$ be a
measuring $\alpha $-pairing. For every left $A$-module $M$ we have:

\begin{enumerate}
\item  $\mathrm{Rat}^{C}(_{A}M)\subseteq M$ is an $A$-submodule.

\item  For every $A$-submodule $N\subset M$ we have $\mathrm{Rat}%
^{C}(_{A}N)=N\cap \mathrm{Rat}^{C}(_{A}M).$

\item  $\mathrm{Rat}^{C}(\mathrm{Rat}^{C}(_{A}M))=\mathrm{Rat}^{C}(_{A}M).$

\item  For every $L\in $ $_{A}\mathcal{M}$ and $f\in $ $\mathrm{Hom}%
_{A-}(M,L)$ we have $f(\mathrm{Rat}^{C}(_{A}M))\subseteq \mathrm{Rat}%
^{C}(_{A}L).$
\end{enumerate}
\end{lemma}

\begin{theorem}
\label{cor-dicht}\emph{(\cite[Lemmata 2.2.8, 2.2.9, Satz 2.2.16]{Abu2001})}
Let $P=(${$A$}$,${$C$}$)$ be a measuring $R$-pairing. Then $\mathcal{M}%
^{C}\subseteq $ $_{A}\mathcal{M}$ and $^{C}\mathcal{M}\subseteq \mathcal{M}%
_{A}$ \emph{(not necessarily full subcategories)}. Moreover the following
are equivalent:

\begin{enumerate}
\item  $P$ satisfies the $\alpha $-condition;

\item  $_{R}C$ is locally projective and $\kappa _{P}(A)\subseteq C^{\ast }$
is dense.

If these equivalent conditions are satisfied, then $\mathcal{M}^{C}\subseteq 
$ $_{A}\mathcal{M}$ and $^{C}\mathcal{M}\subseteq \mathcal{M}_{A}$ are \emph{%
full} subcategories and we have category isomorphisms 
\begin{equation}
\begin{tabular}{lllll}
$\mathcal{M}^{{C}}$ & $\simeq $ & $\mathrm{Rat}^{{C}}(_{A}\mathcal{M})$ & $=$
& $\sigma \lbrack _{{A}}${$C$}$]$ \\ 
& $\simeq $ & $\mathrm{Rat}^{{C}}(_{C^{\ast }}\mathcal{M})$ & $=$ & $\sigma
\lbrack _{C^{\ast }}${$C$}$]$%
\end{tabular}
\text{ and } 
\begin{tabular}{lllll}
$^{{C}}\mathcal{M}$ & $\simeq $ & $^{{C}}\mathrm{Rat}(\mathcal{M}_{A})$ & $=$
& $\sigma \lbrack ${$C$}$_{{A}}]$ \\ 
& $\simeq $ & $^{{C}}\mathrm{Rat}(\mathcal{M}_{C^{\ast }})$ & $=$ & $\sigma
\lbrack ${$C$}$_{C^{\ast }}].$%
\end{tabular}
\label{MC=}
\end{equation}
\end{enumerate}
\end{theorem}

\begin{corollary}
\label{conseq}Let $Q=(B,C)\in \mathcal{P}_{m},$ $\xi :A\rightarrow B$ be an $%
R$-algebra morphism and consider the induced measuring $R$-pairing $%
P:=(A,C). $ Then the following statements are equivalent:

(i) $P\in \mathcal{P}_{m}^{\alpha };$

(ii) $Q\in \mathcal{P}_{m}^{\alpha }$ and $\xi (A)\subset B$ is dense \emph{(%
}w.r.t. the left $C$-adic topology $\mathcal{T}_{C-}(B)$\emph{)};

(iii) $C$ satisfies the $\alpha $-condition and $\kappa _{P}(A)\subset
C^{\ast }$ is dense.

If these equivalent conditions are satisfied, then we get category
isomorphisms 
\begin{equation}
\begin{tabular}{lllll}
$\mathcal{M}^{C}$ & $\simeq $ & $\mathrm{Rat}^{C}(_{A}\mathcal{M})$ & $=$ & $%
\sigma \lbrack _{A}C]$ \\ 
& $\simeq $ & $\mathrm{Rat}^{C}(_{C^{\ast }}\mathcal{M})$ & $=$ & $\sigma
\lbrack _{C^{\ast }}C]$ \\ 
& $\simeq $ & $\mathrm{Rat}^{C}(_{B}\mathcal{M})$ & $=$ & $\sigma \lbrack
_{B}C]$%
\end{tabular}
\text{ and } 
\begin{tabular}{lllll}
$^{C}\mathcal{M}$ & $\simeq $ & $^{C}\mathrm{Rat}(\mathcal{M}_{A})$ & $=$ & $%
\sigma \lbrack C_{A}]$ \\ 
& $\simeq $ & $^{C}\mathrm{Rat}(\mathcal{M}_{C^{\ast }})$ & $=$ & $\sigma
\lbrack C_{C^{\ast }}]$ \\ 
& $\simeq $ & $^{C}\mathrm{Rat}(\mathcal{M}_{B})$ & $=$ & $\sigma \lbrack
C_{B}].$%
\end{tabular}
\label{cat-eq}
\end{equation}
\end{corollary}

\qquad

\begin{punto}
\label{End(C)}Let $(C,\Delta _{C},\varepsilon _{C})$ be an $R$-coalgebra and
denote with $\mathrm{End}^{C}(C)$ (resp. $^{C}\mathrm{End}(C)$) the ring of
all right (resp. left) $C$-colinear morphisms from $C$ to $C$ with the usual
composition. For every right $C$-comodule $M$ we have an isomorphism of $R$%
-modules 
\begin{equation}
\Psi :M^{\ast }\rightarrow \mathrm{Hom}^{C}(M,C),\text{ }h\mapsto \lbrack
m\mapsto \sum f(m_{<0>})m_{<1>}]  \label{M*=End}
\end{equation}
with inverse $g\mapsto \varepsilon _{C}\circ g.$ Analogously $N^{\ast
}\simeq $ $^{C}\mathrm{Hom}(N,C)$ as $R$-modules for every left $C$-comodule 
$N.$ In particular $C^{\ast }\simeq \mathrm{End}^{C}(C)^{op}$ and $C^{\ast
}\simeq $ $^{C}\mathrm{End}(C)$ as $R$-algebras.

If $(A,C)$ is a measuring $\alpha $-pairing, then we have $R$-algebra
isomorphisms 
\begin{equation*}
\mathrm{Biend}(_{A}C):=\mathrm{End}(C_{\mathrm{End}(_{A}C)^{op}})\simeq 
\mathrm{End}(C_{\mathrm{End}^{C}(C)^{op}})\simeq \mathrm{End}(C_{C^{\ast }})=%
\text{ }^{C}\mathrm{End}(C)\simeq C^{\ast }
\end{equation*}
and 
\begin{equation*}
\mathrm{Biend}(C_{A}):=\mathrm{End}(_{\mathrm{End}(C_{A})}C)^{op}\simeq 
\mathrm{End}(_{^{C}\mathrm{End}(C)}C)^{op}\simeq \mathrm{End}(_{C^{\ast
}}C)^{op}=\mathrm{End}^{C}(C)^{op}\simeq C^{\ast }.
\end{equation*}
In particular, if $_{R}C$ is locally projective, then $\mathrm{Biend}%
(_{C^{\ast }}C)\simeq C^{\ast }\simeq \mathrm{Biend}(C_{C^{\ast }})$ as $R$%
-algebras (i.e. $_{C^{\ast }}C_{C^{\ast }}$ is faithfully balanced).
\end{punto}

\begin{corollary}
\label{C=A^r}Let $P=(A,C)\in \mathcal{P}_{m}^{\alpha }$ and consider $%
A^{\ast }$ as an $A$-bimodule with the \emph{regular }$A$-actions 
\begin{equation}
(af)(\widetilde{a})=f(\widetilde{a}a)\text{ and }(fa)(\widetilde{a})=f(a%
\widetilde{a}).  \label{A*-bimodule}
\end{equation}

\begin{enumerate}
\item  For every unitary left \emph{(}right\emph{)} $A$-submodule $%
D\subseteq A^{\ast }$ we have 
\begin{equation*}
\mathrm{Rat}^{C}(_{A}D)=C\cap D\text{\emph{(}}^{C}\mathrm{Rat}(D_{A})=C\cap D%
\text{\emph{)}.}
\end{equation*}
In particular $\mathrm{Rat}^{C}(_{A}A^{\ast })=C=$ $^{C}\mathrm{Rat}%
(A_{A}^{\ast }).$

\item  If $D\subseteq A^{\ast }$ is an $A$-subbimodule, then $_{A}D$ is $C$%
-rational, iff $D_{A}$ is $C$-rational.

\item  Let $R$ be noetherian. If $_{A}A^{\circ }$ \emph{(}equivalently $%
A_{A}^{\circ }$\emph{)} is $C$-rational, then $C=A^{\circ }.$
\end{enumerate}
\end{corollary}

\begin{Beweis}
\begin{enumerate}
\item  Let $D\subset A^{\ast }$ be a left $A$-submodule. By Lemma \ref{clos}
(2) $C\cap D$ is a $C$-rational left $A$-module, i.e. $C\cap D\subseteq 
\mathrm{Rat}^{C}(_{A}D).$ On the other hand, if $f\in \mathrm{Rat}%
^{C}(_{A}D) $ with $\varrho _{D}(f)=\sum f_{i}\otimes c_{i}\in D\otimes
_{R}C,$ then we have for every $a\in A:$%
\begin{equation*}
f(a)=(af)(1_{A})=\sum f_{i}(1_{A})<a,c_{i}>,
\end{equation*}
i.e. $f=\sum f_{i}(1_{A})c_{i}\in C.$ Hence $\mathrm{Rat}^{C}(_{A}D)=C\cap
D. $ The corresponding result for right $A$-submodules $D\subset A^{\ast }$
follows by symmetry.

\item  Let $D\subseteq A^{\ast }$ be an $A$-subbimodule. Then by (1) $%
\mathrm{Rat}^{C}(_{A}D)=C\cap D=$ $^{C}\mathrm{Rat}(D_{A}).$

\item  If $R$ is noetherian, then $A^{\circ }\subset A^{\ast }$ is an $A$%
-subbimodule. Obviously $C\overset{\chi _{P}}{\hookrightarrow }A^{\circ }$
and it follows by assumption and (1) that $A^{\circ }=\mathrm{Rat}%
^{C}(_{A}A^{\circ })=C\cap A^{\circ }=C.\blacksquare $
\end{enumerate}
\end{Beweis}

An important rule by the study of the category of rational representations
of measuring $\alpha $-pairings is played by the

\begin{punto}
\textbf{Finiteness Theorem.}\label{es}

\begin{enumerate}
\item  Let $P=(A,C)$ be a measuring $\alpha $-pairing. If $M\in \mathrm{Rat}%
^{C}(_{A}\mathcal{M}$) (resp. $M\in $ $^{C}\mathrm{Rat}(\mathcal{M}_{A}),$ $%
M\in $ $^{C}\mathrm{Rat}^{C}(_{A}\mathcal{M}_{A})$), then there exists for
every finite set $\{m_{1},...,m_{k}\}\subset M$ some $N\in \mathrm{Rat}%
^{C}(_{A}\mathcal{M}$) (resp. $N\in $ $^{C}\mathrm{Rat}(\mathcal{M}_{A}),$ $%
N\in $ $^{C}\mathrm{Rat}^{C}(_{A}\mathcal{M}_{A})$), such that $N_{R}$ is
finitely generated.

\item  Let $C$ be a locally projective $R$-coalgebra. Then every finite
subset of $C$ is contained in a right $C$-coideal (resp. a left $C$-coideal,
a $C$-bicoideal), that is finitely generated in $\mathcal{M}_{R}.$
\end{enumerate}
\end{punto}

\begin{Beweis}
\begin{enumerate}
\item  Let $P=(A,C)\in \mathcal{P}_{m}^{\alpha }.$ Let $M\in \mathrm{Rat}%
^{C}(_{A}\mathcal{M}$) and $\{m_{1},...,m_{k}\}\subset M.$ Then $%
Am_{i}\subset M$ is an $A$-submodule, hence a $C$-subcomodule. Moreover $%
m_{i}\in Am_{i}$ and so there exists a subset $\{(m_{ij},c_{ij})%
\}_{j=1}^{n_{i}}\subset Am_{i}\times C,$ such that $\varrho
_{M}(m_{i})=\sum\limits_{j=1}^{n_{i}}m_{ij}\otimes c_{ij}$ for $i=1,...,k.$
Obviously $N:=\sum\limits_{i=1}^{k}Am_{i}=\sum\limits_{i=1}^{k}\sum%
\limits_{j=1}^{n_{i}}Rm_{ij}\subset M$ is a $C$-subcomodule and contains $%
\{m_{1},...,m_{k}\}.$

Using analogous arguments one can show the corresponding result for $C$%
-rational right $A$-modules and $C$-birational $A$-bimodules.

\item  If $C$ is a locally projective $R$-coalgebra, then $(C^{\ast },C)\in 
\mathcal{P}_{m}^{\alpha }$ and the result follows by (1).$\blacksquare $
\end{enumerate}
\end{Beweis}

\qquad The following result gives topological characterizations of the $C$%
-rational left $A$-modules and generalizes the corresponding result obtained
by D. Radford \cite[2.2]{Rad73} from the case of base fields to the case of
arbitrary (artinian) commutative ground rings (see also \cite[Proposition
1.4.4]{LR97}).

\begin{proposition}
\label{n_rat}Let $P=(A,C)$ be a measuring $\alpha $-pairing and consider $A$
with the $C$-adic topology $\mathcal{T}_{C-}(A)=A[\frak{T}_{ls}(C)].$ If $M$
is a unitary left $A$-module, then for every $m\in M$ the following
statements are equivalent:

\begin{enumerate}
\item  there exists a finite subset $W=\{c_{1},...,c_{k}\}\subset C$, such
that $(0_{C}:W)\subset (0_{M}:m).$

\item  $Am$ is $C$-subgenerated;

\item  $m\in \mathrm{Rat}^{C}(_{A}M).$

\item  there exists a finitely generated $R$-submodule $K\subset C,$ such
that $K^{\bot }\subset (0_{M}:m).$

If $R$ is artinian, then ``1-4'' are moreover equivalent to:

\item  $(0_{M}:Am)$ contains an $R$-cofinite closed $R$-submodule of $A;$

\item  $(0_{M}:Am)$ is an $R$-cofinite closed $A$-ideal;

\item  $(0_{M}:m)$ contains an $R$-cofinite closed $A$-ideal;

\item  $(0_{M}:m)$ is an $R$-cofinite closed left $A$-ideal.
\end{enumerate}
\end{proposition}

\begin{Beweis}
(1) $\Rightarrow $ (2) By assumption and \ref{C-ad} $m\in N:=\mathrm{Sp}%
(\sigma \lbrack _{A}C],M).$ Since $Am\subset $ $N$ is an $A$-submodule, it's 
$C$-subgenerated.

(2) $\Rightarrow $ (3) By assumption and Theorem \ref{cor-dicht} $m\in
Am\subset \mathrm{Rat}^{C}(_{A}M).$

(3)\ $\Rightarrow $ (4) Let $\varrho (m)=\sum\limits_{i=1}^{k}m_{i}\otimes
c_{i}$ and $K:=\sum\limits_{i=1}^{k}Rc_{i}\subset C.$ Then obviously $%
K^{\bot }\subset (0_{M}:m).$

(4) $\Rightarrow $ (1) For every subset $W\subset C$ we have $%
(0_{C}:W)\subseteq W^{\bot }.$

Let $R$ be \emph{artinian}.

(3) $\Rightarrow $ (5). By Theorem \ref{cor-dicht} $\mathrm{Rat}^{C}(_{A}M)$
is a $C$-rational left $A$-module. Assume that $\varrho
_{M}(m)=\sum\limits_{i=1}^{k}m_{i}\otimes c_{i}\in \mathrm{Rat}%
^{C}(_{A}M)\otimes _{R}C,$ $\varrho
_{M}(m_{i})=\sum\limits_{j=1}^{n_{i}}m_{ij}\otimes c_{ij}$ for $i=1,...,k$
and set $K:=\sum\limits_{i=1}^{k}\sum\limits_{j=1}^{n_{i}}Rc_{ij}.$ Then we
have for every $a\in K^{\bot }$ and arbitrary $b\in A:$%
\begin{equation*}
a(bm)=a(\sum\limits_{i=1}^{k}m_{i}<b,c_{i}>)=\sum\limits_{i=1}^{k}\sum%
\limits_{j=1}^{n_{i}}m_{ij}<a,c_{ij}><b,c_{i}>=0,
\end{equation*}
i.e. $K^{\bot }\subseteq (0_{M}:Am).$ The $R$-module $K$ is finitely
generated and it follows from the embedding $A/K^{\bot }\hookrightarrow
K^{\ast },$ that $K^{\bot }\subset A$ is an $R$-cofinite $R$-submodule.
Moreover $K^{\bot }$ is by \cite[Lemma 1.7 (1)]{Abu} closed.

If $R$ is artinian, then the implications (5) $\Rightarrow $ (6) $%
\Rightarrow $ (7) $\Rightarrow $ (8) $\Rightarrow $ (4) follow from 
\cite[Lemma 1.7 (4)]{Abu}.$\blacksquare $
\end{Beweis}

\begin{lemma}
\label{R_f}Let $C$ be an $R$-coalgebra and consider $C^{\ast }$ with the
finite topology. For every $f\in C^{\ast }$ the $R$-linear mappings 
\begin{equation*}
\xi _{f}^{r}:C^{\ast }\rightarrow C^{\ast },\text{ }g\mapsto g\star f\text{
and }\xi _{f}^{l}:C^{\ast }\rightarrow C^{\ast },\text{ }g\mapsto f\star g
\end{equation*}
are continuous. If $R$ is an injective cogenerator, then $\xi _{f}^{r}$ and $%
\xi _{f}^{l}$ are linearly closed \emph{(}i.e. $\xi _{f}^{r}(X)\subset
C^{\ast }$ and $\xi _{f}^{l}(X)\subset C^{\ast }$ are closed for every
closed $R$-submodule $X\subset C^{\ast }$\emph{)}.
\end{lemma}

\begin{Beweis}
Consider for every $f\in C^{\ast }$ the $R$-linear mappings 
\begin{equation*}
\theta _{f}^{l}:C\rightarrow C,c\mapsto f\rightharpoonup c\text{ and }\theta
_{f}^{r}:C\rightarrow C,\text{ }c\mapsto c\leftharpoonup f\text{.}
\end{equation*}
Then we have for every $g\in C^{\ast }$ and $c\in C:$%
\begin{equation*}
\xi _{f}^{r}(g)(c)=(g\star f)(c)=\sum g(c_{1})f(c_{2})=g(f\rightharpoonup
c)=g(\theta _{f}^{l}(c))=((\theta _{f}^{l})^{\ast }(g))(c).
\end{equation*}
So $\xi _{f}^{r}=(\theta _{f}^{l})^{\ast }$ and analogously $\xi
_{f}^{l}=(\theta _{f}^{r})^{\ast }.$ The result follows then by 
\cite[Proposition 1.10]{Abu}.$\blacksquare $
\end{Beweis}

\qquad If $P=(A,C)\in \mathcal{P}_{m}^{\alpha }$, then the Grothendieck
category $\mathrm{Rat}^{C}(_{A}\mathcal{M})\simeq \sigma \lbrack _{A}C]$ is
in general not closed under extensions:

\begin{example}
(\cite[Page 520]{Rad73}) Let $R$ be a \emph{base field}, $V$ be an infinite
dimensional vector space over $R$ and consider the $R$-coalgebra $C:=R\oplus
V$ (with $\Delta (v)=1\otimes v+v\otimes 1$ and $\varepsilon (v)=1_{R}$).
Let $I\subset V^{\ast }$ be vector subspace that is not closed, and consider
the exact sequence of $C^{\ast }$-modules 
\begin{equation*}
0\rightarrow V^{\ast }/I\rightarrow C^{\ast }/I\rightarrow C^{\ast }/V^{\ast
}\rightarrow 0.
\end{equation*}
Then $V^{\ast }/I$ and $C^{\ast }/V^{\ast }$ are $C$-rational, while $%
C^{\ast }/I$ is not.
\end{example}

\begin{lemma}
\label{jm}\emph{(\cite[Lemma 6.1.1, Corollary 6.1.2]{Swe69})} Let $%
I\vartriangleleft A$ be an ideal.

\begin{enumerate}
\item  Let $M$ be a finitely generated left \emph{(}right\emph{)} $A$%
-module. If $_{A}I$ \emph{(}$I_{A}$\emph{)} is finitely generated, then also 
$IM\subset M$ \emph{(}$MI\subset M$\emph{)} is a finitely generated $A$%
-submodule. If $I\subset A$ is $R$-cofinite, then $IM\subset M$ \emph{(}$%
MI\subset M$\emph{)} is an $R$-cofinite $A$-submodule.

\item  If $_{A}I$ \emph{(}$I_{A}$\emph{)} is finitely generated, then $%
_{A}I^{n}$ \emph{(}$I_{A}^{n}$\emph{)} is finitely generated for every $%
n\geq 1.$ If moreover $I\subset A$ is $R$-cofinite, then $I^{n}\subset A$ is 
$R$-cofinite.
\end{enumerate}
\end{lemma}

\qquad The following result generalizes \cite[2.5]{Rad73} from the case of
base fields to the case of arbitrary commutative QF\ rings.

\begin{proposition}
\label{part}Let $R$ be a QF Ring, $C$ be a projective $R$-coalgebra and
consider an exact sequence of left $C^{\ast }$-modules 
\begin{equation*}
0\longrightarrow N\overset{\iota }{\longrightarrow }M\overset{\pi }{%
\longrightarrow }L\longrightarrow 0.
\end{equation*}
If $N,L\in \mathrm{Rat}^{C}(_{C^{\ast }}\mathcal{M})$ and $_{C^{\ast }}(0:l)$
is finitely generated for every $l\in L,$ then $M$ is $C$-rational.
\end{proposition}

\begin{Beweis}
Let $m\in M$ and $\{f_{1},...,f_{k}\}$ be a generating system of $_{C^{\ast
}}(0_{L}:\pi (m)).$ By assumption $\pi (m)$ is $C$-rational and so there
exist by Proposition \ref{n_rat} $R$-cofinite closed $A$-ideals $%
J_{i}\subset (0_{N}:f_{i}m)$ for $i=1,...,k.$ So we have for the closed $R$%
-cofinite $A$-ideal $J:=\bigcap\limits_{i=1}^{k}J_{i}\vartriangleleft
C^{\ast }:$%
\begin{equation*}
J(0_{L}:\pi (m))\rightharpoonup m=(J\star f_{1}+...+J\star
f_{k})\rightharpoonup m=0,
\end{equation*}
i.e. $J(0_{L}:\pi (m))\subseteq (0_{M}:m).$ By Lemmata \ref{R_f}, \ref{jm}
and \cite[Proposition 1.10 (3.d)]{Abu} $J(0_{L}:\pi
(m))=\sum\limits_{i=1}^{k}J\star f_{i}$ is $R$-cofinite and closed. It
follows then by \cite[Lemma 1.7 (4)]{Abu} that $(0_{M}:m)\vartriangleleft
_{l}C^{\ast }$ is $R$-cofinite and closed, hence $m\in \mathrm{Rat}%
^{C}(_{C^{\ast }}M)$ by Proposition \ref{n_rat}.$\blacksquare $
\end{Beweis}

\begin{definition}
An $R$-algebra $A$ is called \emph{nearly left noetherian} (resp. \emph{%
nearly right noetherian}, \emph{nearly noetherian}), if every $R$-cofinite%
{\normalsize \ }left (resp. right, two-sided) $A$-ideal is finitely
generated in $_{A}\mathcal{M}$ (resp. in $\mathcal{M}_{A},$ in $_{A}\mathcal{%
M}_{A}$).
\end{definition}

\qquad As a corollary of Theorem \ref{cor-dicht} and Proposition \ref{part}
we get

\begin{corollary}
Let $R$ be a QF Ring and $C$ be a projective $R$-coalgebra. If $C^{\ast }$
is nearly left noetherian \emph{(}resp. nearly right noetherian\emph{)},
then $\mathcal{M}^{C}\simeq \mathrm{R}\mathrm{at}^{C}(_{C^{\ast }}C)=\sigma
\lbrack _{C^{\ast }}C]$ \emph{(}resp. $^{C}\mathcal{M}\simeq $ $^{C}\mathrm{%
Rat}(\mathcal{M}_{C^{\ast }})=\sigma \lbrack C_{C^{\ast }}]$\emph{)} is
closed under extensions.
\end{corollary}

\subsection*{Duality relations between substructures}

\qquad As an application of our results in this section and our observations
about the linear weak topology \cite{Abu} we generalize known results on the
duality relations between substructures of a coalgebra and substructures of
its dual algebra from the case of base fields (e.g. \cite{Swe69}, \cite
{Abe80} and \cite[1.5.29]{DNR}) to the case of measuring $\alpha $-pairings
over arbitrary commutative rings.

\qquad As a consequence of Theorem \ref{cor-dicht} and \cite[Theorem 1.8]
{Abu} we get

\begin{proposition}
\label{id-coid}Let $P=(A,C)\in \mathcal{P}_{m}.$

\begin{enumerate}
\item  Let $K\subset C$ be an $R$-submodule.

If $K$ is a right \emph{(}a left\emph{)} $C$-coideal, then $K^{\bot }=%
\mathrm{An}_{A}^{r}(K)$ is a right \emph{(}a left\emph{)} $A$-ideal;

If $K$ is a $C$-bicoideal, then $K^{\bot }=\mathrm{An}_{A}^{r}(K)\cap 
\mathrm{An}_{A}^{l}(K)$ is a two-sided $A$-ideal.

\item  Let $P\in \mathcal{P}_{m}^{\alpha }.$

\begin{enumerate}
\item  For every $R$-submodule $I\subset A$ we have:

If $I\subset A$ is a right \emph{(}a left\emph{)} ideal, then $I^{\bot
}\subset C$ is a right \emph{(}a left\emph{)} coideal;

If $I\vartriangleleft A$ is a two-sided ideal \emph{(}and $I^{\bot }\subset
C $ is pure\emph{)}, then $I^{\bot }\subset C$ is a bicoideal \emph{(}an $R$%
-subcoalgebra\emph{)}.

\item  Let $R$ be an injective cogenerator. For a closed $R$-submodule $%
I\subset A$ we have:

$I$ is a right \emph{(}a left\emph{)} ideal, iff $I^{\bot }\subset C$ is a
right \emph{(}a left\emph{) }coideal.

$I$ is a two-sided ideal \emph{(}and $I^{\bot }\subset C$ is pure\emph{)},
iff $I^{\bot }\subset C$ is a bicoideal \emph{(}an $R$-subcoalgebra\emph{)}.
\end{enumerate}
\end{enumerate}
\end{proposition}

\begin{lemma}
\label{An(B)-coid}

\begin{enumerate}
\item  If $P=(A,C)$ is a measuring $R$-pairing and $K\subset C$ is a
coideal, then $K^{\bot }\subset A$ is an $R$-subalgebra with unity $1_{A}.$

\item  Let $R$ be a QF Ring, $C$ be a projective $R$-coalgebra and $A\subset
C^{\ast }$ be an $R$-subalgebra \emph{(}with $\varepsilon _{C}\in A$\emph{)}%
. If $\mathrm{Ke}(A)\subset C$ is pure, then $\Delta _{C}(\mathrm{Ke}%
(A))\subset \mathrm{Ke}(A)\otimes _{R}C+C\otimes _{R}\mathrm{Ke}(A)$ \emph{(}%
$\mathrm{Ke}(A)\subset C$ is a $C$-coideal\emph{)}.
\end{enumerate}
\end{lemma}

\begin{Beweis}
\begin{enumerate}
\item  Obvious.

\item  Let $A\subset C^{\ast }$ be an $R$-subalgebra and consider the
canonical $R$-linear mappings 
\begin{equation*}
\kappa :A\otimes _{R}A\rightarrow (C\otimes _{R}C)^{\ast }\text{ and }\chi
:C\otimes _{R}C\rightarrow (A\otimes _{R}A)^{\ast }.
\end{equation*}
If $\mathrm{Ke}(A)\subset C$ is pure, then it follows form \cite[Proposition
1.10 (3.c), Corollary 2.9]{Abu} that 
\begin{equation}
\begin{tabular}{lllll}
$\mathrm{Ke}(A)$ & $=$ & $\mathrm{Ke}(\Delta _{C}^{\ast }(\kappa (A\otimes
_{R}A)))$ & $=$ & $\Delta _{C}^{-1}(\mathrm{Ke}(\kappa (A\otimes _{R}A)))$
\\ 
& $=$ & $\Delta _{C}^{-1}(\mathrm{Ke}(A)\otimes _{R}C+C\otimes _{R}\mathrm{Ke%
}(A)),$ &  & 
\end{tabular}
\label{Bij-cood-al}
\end{equation}
i.e. $\Delta _{C}(\mathrm{Ke}(A))\subseteq \mathrm{Ke}(A)\otimes
_{R}C+C\otimes _{R}\mathrm{Ke}(A).$ If moreover $\varepsilon _{C}\in A,$
then $\varepsilon _{C}(\mathrm{Ke}(A))=0,$ i.e. $\mathrm{Ke}(A)\subset C$ is
a $C$-coideal.$\blacksquare $
\end{enumerate}
\end{Beweis}

\qquad As a consequence of Propositions \ref{id-coid}, \ref{An(B)-coid} and 
\cite[Theorem 1.8]{Abu} we get

\begin{corollary}
\label{op-id*}Let $R$ be an injective cogenerator and $C$ a locally
projective $R$-coalgebra. If we denote with $\mathcal{C}\ $the class of all $%
R$-submodules of $C$ and with $\mathcal{H}\ $the class of all $R$-submodules
of $C^{\ast },$ then 
\begin{equation}
\mathrm{An}(-):\mathcal{C}\rightarrow \mathcal{H}\text{{\normalsize \ }and }%
\mathrm{Ke}(-):\mathcal{H}\rightarrow \mathcal{C}  \label{An(-)}
\end{equation}
induce bijections 
\begin{equation}
\begin{tabular}{lll}
$\{K\subset C$ a right $C$-coideal$\}$ & $\leftrightarrow $ & $%
\{I\vartriangleleft _{r}C^{\ast }$ a \emph{closed }right $A$-ideal$\}$ \\ 
$\{K\subset C$ a left $C$-coideal$\}$ & $\leftrightarrow $ & $%
\{I\vartriangleleft _{l}C^{\ast }$ a \emph{closed} left $A$-ideal$\}$ \\ 
$\{K\subset C$ a $C$-bicoideal$\}$ & $\leftrightarrow $ & $%
\{I\vartriangleleft C^{\ast }\text{ a \emph{closed}}$\emph{\ }two-sided ideal%
$\}.$ \\ 
$\{K\subset C\text{ an }R\text{-subcoalgebra}\}$ & $\leftrightarrow $ & $%
\{I\vartriangleleft C^{\ast }$ a \emph{closed\ }two-sided ideal, $\mathrm{Ke}%
(I)\subset C\text{ pure}\}.$%
\end{tabular}
\label{C-C*}
\end{equation}
If $R$ is moreover a QF ring, then \emph{(\ref{An(-)})} induces a bijection 
\begin{equation*}
\{K\subset C\text{ a pure }C\text{-coideal}\}\leftrightarrow \{A\subset
C^{\ast }\text{ a }\emph{closed\ }R\text{-subalgebra, }\varepsilon _{C}\in A,%
\text{ }\mathrm{Ke}(A)\subset C\text{ pure}\}.
\end{equation*}
\end{corollary}

\section{Dual coalgebras}

Every $R$-coalgebra $(C,\Delta _{C},\varepsilon _{C})$ has a \emph{dual }$R$%
\emph{-algebra,} namely $C^{\ast }$ with multiplication the convolution
product 
\begin{equation*}
\star :C^{\ast }\otimes _{R}C^{\ast }\overset{\delta }{\longrightarrow }%
(C\otimes _{R}C)^{\ast }\overset{\Delta _{C}^{\ast }}{\longrightarrow }%
C^{\ast },
\end{equation*}
where $\delta $ is the canonical $R$-linear mapping, and with unity element $%
\varepsilon _{C}.$ If $(A,\mu _{A},\eta _{A})$ is an $R$-algebra that is
finitely generated projective as an $R$-module, then $A^{\ast }$ becomes an $%
R$-coalgebra with comultiplication given by 
\begin{equation*}
\mu _{A}^{\circ }:A^{\ast }\overset{\mu _{A}^{\ast }}{\longrightarrow }%
(A\otimes _{R}A)^{\ast }\overset{\delta ^{-1}}{\longrightarrow }A^{\ast
}\otimes _{R}A^{\ast },
\end{equation*}
where $\delta :A^{\ast }\otimes _{R}A^{\ast }\rightarrow (A\otimes
_{R}A)^{\ast }$ is the canonical isomorphism, and with counity $\eta
_{A}^{\ast }:A^{\ast }\rightarrow R.$ If $A$ is not finitely generated
projective, then $\delta $ in not surjective anymore (and not even injective
over arbitrary ground rings), hence $\mu _{A}^{\circ }$ is not well defined
and $\mu _{A}$ incudes on $A^{\ast }$ no $R$-coalgebra structure. However,
if $R$ is base field and we consider the $R$-algebra $A$ with the cofinite
topology $\mathrm{Cf}(A)$ (see \ref{ko-top}) and $R$ with the discrete
topology, then the character module\ $A^{\circ }$ of all continuous $R$%
-linear mappings from $A$ to $R$ is an $R$-coalgebra (\cite{Swe69}). That
result was generalized in \cite{CN90} to the case of \emph{Dedekind domains}
and in \cite{AG-TW2000} to the case of arbitrary \emph{noetherian }(\emph{%
hereditary}) commutative rings.

In this section we consider coalgebra structures on the character module of
an algebra, considered with a linear topology induced from a filter basis
consisting of cofinite ideals over an arbitrary (noetherian) ring.

\begin{punto}
Let $A$ be an $R$-algebra and $\frak{B}$ be a filter basis consisting of $R$%
-cofinite $A$-ideals. Then $\frak{B}$ induces on $A$ a left linear topology $%
\frak{T}(\frak{B}),$ such that $(A,\frak{T}(\frak{B}))$ is a left linear
topological $R$-algebra and $\frak{B}\ $is a neighbourhood basis of $0_{A}.$
With 
\begin{equation}
A_{\frak{B}}^{\circ }:=\{f\in A^{\ast }|\text{ }\exists \text{ }I\in \frak{B}%
,\text{ such that }f(I)=0\}={\underrightarrow{lim}}_{\frak{B}}(A/I)^{\ast }
\label{A_F-0}
\end{equation}
we denote the \emph{character module} of all continuous $R$-linear mappings
from $A$ to $R$ (where $R$ is considered as usual with the discrete
topology). The \emph{completion} of $A$ w.r.t. $\frak{B}$ is denoted with 
\begin{equation*}
\widehat{A}_{\frak{B}}:=\underleftarrow{lim}\{A/I|\text{ }I\in \frak{B}\}.
\end{equation*}
If $A_{\frak{B}}^{\circ }$ is an $R$-coalgebra, then we call $A_{\frak{B}%
}^{\circ }$ the \emph{continuous dual }$R$\emph{-coalgebra of }$A$\emph{\
w.r.t.}\textbf{\ }$\frak{B}.$

Analogously $\frak{B}$ induces on $A$ a right linear topology, such that $A$
is a right linear topological $R$-algebra and $\frak{B}$ is a neighbourhood
basis of $0_{A}.$
\end{punto}

\begin{remark}
\label{twoside}(Compare \cite[Proposition 3.1]{CG-RTvO2001}) Let $R$ be
noetherian and $A$ be an $R$-algebra. Let $I$ be an $R$-cofinite left $A$%
-ideal, say $A/I=\sum\limits_{i=1}^{k}R(a_{i}+I),$ and consider the \emph{%
two-sided} $A$-ideal 
\begin{equation*}
J:=\bigcap_{i=1}^{k}(I:a_{i})=(I:A)\subset (I:1_{A})=I.
\end{equation*}
Then 
\begin{equation*}
\varphi _{I}:A\rightarrow \mathrm{End}_{R}(A/I),\text{ }a\mapsto \lbrack
b+I\mapsto ab+I]
\end{equation*}
is an $R$-algebra morphism with $\mathrm{Ke}(\varphi _{I})=J,$ i.e. $J$ is
an $R$-cofinite $A$-ideal.

Analogously one can show that every $R$-cofinite right $A$-ideal contains an 
$R$-cofinite \emph{two-sided} $A$-ideal.$\blacksquare $
\end{remark}

\qquad The following result extends \cite[1.11]{AG-TW2000} and \cite[Remark
2.14]{AG-TL2001}:

\begin{theorem}
\label{R(G)-co}Let $R$ be noetherian and $A$ an $R$-algebra. If $C\subseteq
A^{\circ }$ is an $A$\emph{-subbimodule} and $P:=(A,C),$ then the following
statements are equivalent:

\begin{enumerate}
\item  $_{R}C$ is locally projective and $\kappa _{P}(A)\subset C^{\ast }$
id dense.

\item  $_{R}C$ satisfies the $\alpha $-condition and $\kappa _{P}(A)\subset
C^{\ast }$ is dense;

\item  $(A,C)$ is an $\alpha $-pairing;

\item  $C\subset R^{A}$ is pure \emph{(}in the sense of Cohn\emph{)};

\item  $C$ is an $R$-coalgebra and $(A,C)\in \mathcal{P}_{\alpha }^{m}.$

If $R$ is a QF Ring, then ``1-4'' are moreover equivalent to

\item  $_{R}C$ is projective.
\end{enumerate}
\end{theorem}

\begin{Beweis}
The equivalences (1) $\Leftrightarrow $ (2) and (3) $\Leftrightarrow $ (4)
follow from \cite[Lemma 2.13, Proposition 2.6 (3), ]{Abu}.

(2)\ $\Rightarrow $ (3) follows from \cite[Proposition 2.4 (2) ]{Abu}.

(4) $\Rightarrow $ (5) If $C\subset R^{A}$ is pure, then $C$ is by 
\cite[1.11]{AG-TW2000} an $R$-coalgebra. It follows moreover for all $f\in C$
and arbitrary $a,\widetilde{a}\in A$ that 
\begin{equation*}
\kappa _{P}(a\widetilde{a})(f)=f(a\widetilde{a})=\sum f_{1}(a)f_{2}(%
\widetilde{a})=(\kappa _{P}(a)\underline{\otimes }\kappa _{P}(\widetilde{a}%
))(\Delta (f))=(\kappa _{P}(a)\star \kappa _{P}(\widetilde{a}))(f)
\end{equation*}
and 
\begin{equation*}
\kappa _{P}(1_{A})(f)=f(1_{A})=\varepsilon _{C}(f)\text{ for all }f\in C.
\end{equation*}
So $\kappa _{P}:A\rightarrow C^{\ast }$ is an $R$-algebra morphism, i.e. $%
P\in \mathcal{P}_{m}.$ By \cite[Proposition 2.6]{Abu} $P$ satisfies the $%
\alpha $-condition, hence $P\in \mathcal{P}_{\alpha }^{m}.$

(5) $\Rightarrow $ (2) follows from Theorem \ref{cor-dicht}.

Let $R$ be a QF ring.

(2) $\Rightarrow $ (6) follows from Remark \ref{flat}.

(6) $\Rightarrow $ (2) If $_{R}C$ is projective, then $C\;$satisfies the $%
\alpha $-condition by \cite[Proposition 2.14 (3)]{Abu}. Consider the $R$%
-submodule $\kappa _{P}(A)\subset C^{\ast }.$ By \cite[Theorem 1.8 (1)]{Abu}
we have 
\begin{equation*}
\overline{\kappa _{P}(A)}:=\mathrm{An}{\mathrm{Ke}}(\kappa _{P}(A))=\mathrm{%
An}(A^{\bot })=\mathrm{An}(0_{C})=C^{\ast },
\end{equation*}
i.e. $\kappa _{P}(A)\subset C^{\ast }$ is dense.$\blacksquare $
\end{Beweis}

\begin{definition}
An $R$-algebra $A$ is said to \emph{satisfy the }$\alpha $\emph{-condition, }%
if the class $\mathcal{K}_{A}$ of all $R$-cofinite $A$-ideals is a filter
and the induced $R$-pairing $(A,A^{\circ })$ satisfies the $\alpha $%
-condition (in case $R$ is noetherian, this is equivalent to the purity of $%
A^{\circ }\subset R^{A}$). An $R$-coalgebra $C$ is said to \emph{satisfy the 
}$\alpha $\emph{-condition} or to be an $\alpha $\emph{-coalgebra}, if the $%
R $-pairing $(C^{\ast },C)$ satisfies the $\alpha $-condition (equivalently,
if $_{R}C$ is locally projective). With $\mathbf{Alg}_{R}^{\alpha }\subseteq 
\mathbf{Alg}_{R}$ resp. $\mathbf{Cog}_{R}^{\alpha }\subseteq \mathbf{Cog}%
_{R} $ we denote the \emph{full }subcategory of $\alpha $-algebras resp. $%
\alpha $-coalgebras.
\end{definition}

\begin{remark}
\label{C_A0} Let $R$ be noetherian and $A$ be an $\alpha $-algebra. Then
there is obviously a 1-1 correspondence 
\begin{equation*}
\{P=(A,C)|\text{ }P\in \mathcal{P}{\normalsize _{m}^{\alpha }}%
\}\longleftrightarrow \{C|\text{ }C\subseteq A^{\circ }\text{ is an }R\text{%
-subcoalgebra}\}.
\end{equation*}
\end{remark}

\begin{lemma}
\label{alg-coalg}

\begin{enumerate}
\item  If $C,$ $D$ are $R$-coalgebras and $\theta :D\rightarrow C$ is an $R$%
-coalgebra morphism, then $\theta ^{\ast }:C^{\ast }\rightarrow D^{\ast }$
is an $R$-algebra morphism and 
\begin{equation*}
(\theta ^{\ast },\theta ):(D^{\ast },D)\rightarrow (C^{\ast },C)
\end{equation*}
is a morphism in $\mathcal{P}_{m}.$

\item  Let $R$ be noetherian, $A,B$ be $\alpha $-algebras and $\xi
:A\rightarrow B$ be an $R$-algebra morphism. Then we have a morphism in $%
\mathcal{P}_{m}^{\alpha }$ 
\begin{equation*}
(\xi ,\xi ^{\circ }):(B,B^{\circ })\rightarrow (A,A^{\circ }).
\end{equation*}
\end{enumerate}
\end{lemma}

\begin{Beweis}
\begin{enumerate}
\item  Trivial.

\item  If $f\in B^{\circ },$ then there exists an $R$-cofinite $B$-ideal $%
I\vartriangleleft B,$ such that $f\in (B/I)^{\ast }.$ By assumption $R$ is
noetherian and so $\xi ^{-1}(I)\subset A$ is an $R$-cofinite $A$-ideal, i.e. 
$\xi ^{\circ }(f)\in A^{\circ }$ and we get a morphism of $R$-pairings 
\begin{equation*}
(\xi ,\xi ^{\circ }):(B,B^{\circ })\rightarrow (A,A^{\circ }).
\end{equation*}
By assumption $\xi $ is an $R$-algebra morphism. Moreover the canonical $R$%
-linear mapping $A^{\circ }\otimes _{R}A^{\circ }\rightarrow (A\otimes
_{R}A)^{\ast }$ is by \cite[Corollary 2.8 (1)]{Abu} an embedding, hence $\xi
^{\circ }:B^{\circ }\rightarrow A^{\circ }$ is an $R$-coalgebra morphism by
Lemma \ref{prop-mes} (1).$\blacksquare $
\end{enumerate}
\end{Beweis}

\begin{lemma}
Let $R$ be noetherian, $B$ an $\alpha $-algebra and consider the $\alpha $%
-pairing $(B,B^{\circ }).$ If $A\subset B$ is an $\alpha $-subalgebra with $%
1_{B}\in A,$ then $A^{\bot }:=\mathrm{An}(A)\cap B^{\circ }$ is a $B^{\circ
} $-coideal.
\end{lemma}

\begin{Beweis}
The embedding $\iota _{A}:A\hookrightarrow B$ is an $R$-algebra morphism and
so $\iota _{A}^{\circ }:B^{\circ }\rightarrow A^{\circ }$ is by Lemma \ref
{prop-mes} (1) an $R$-coalgebra morphism. Hence $A^{\bot }:=\mathrm{Ke}({%
\iota _{A}^{\circ }})\subset B^{\circ }$ is a $B^{\circ }$-coideal.$%
\blacksquare $
\end{Beweis}

\qquad The following result follows directly from Propositions \ref{id-coid}%
, \ref{An(B)-coid}, Lemma \ref{alg-coalg} and \cite[Theorem 1.8]{Abu}:

\begin{corollary}
Let $R$ be a QF Ring, $A$ be an $\alpha $-algebra, $P:=(A,A^{\circ })$ and
consider $A$ with the linear weak topology $A[\frak{T}_{ls}(A^{\circ })].$%
{\normalsize \ }Let $I\subset A$ be a closed $R$-submodule and set $I^{\bot
}:=\mathrm{An}(I)\cap A^{\circ }.$ Then $I$ is a right \emph{(}a left\emph{)}
$A$-ideal, iff $I^{\bot }$ is a right \emph{(}a left\emph{) }$A^{\circ }$%
-coideal. Moreover $I\subset A$ is a two-sided $A$-ideal \emph{(}and $%
I^{\bot }\subset A^{\circ }$ is pure\emph{)}, iff $I^{\bot }\subset A^{\circ
}$ is an $A^{\circ }$-bicoideal \emph{(}an $R$-subcoalgebra\emph{)}.
\end{corollary}

\subsection*{The convolution coalgebra}

Dual to the convolution algebra, D. Radford presented in \cite{Rad73} the so
called \emph{convolution coalgebra} in the case of base fields. Over
arbitrary noetherian ground rings the following version of his definition
makes sense:

\begin{punto}
\label{konv-ko} Let $R$ be noetherian. If $C$ is an $R$-coalgebra and $A$ is
an $\alpha $-algebra, then we call $A\star C:=A^{\circ }\otimes _{R}C$ the 
\emph{convolution coalgebra} of $A$ and $C.$ In the special case $C=R$ we
have $A\star R\simeq A^{\circ }.$
\end{punto}

\qquad The following result generalizes results of D. Radford \cite{Rad73}
on the convolution coalgebra from the case of base fields to the case of
arbitrary noetherian ground rings:

\begin{punto}
\label{con-con}Let $R$ be noetherian, $C$ be a locally projective $R$%
-coalgebra and $A$ be an $\alpha $-algebra. It's easy to see then that $%
P:=(A\otimes _{R}C^{\ast },A\star C)$ is a measuring $R$-pairing, which
satisfies the $\alpha $-condition by \cite[Lemma 2.8]{Abu}. By \cite{Rad73}
the following mappings are $R$-algebra morphisms: 
\begin{equation*}
\begin{tabular}{llllllll}
$\beta $ & $:$ & $\mathrm{Hom}_{R}(C,A)$ & $\rightarrow $ & $(A^{\circ
}\otimes _{R}C)^{\ast },$ & $f$ & $\mapsto $ & $[h\otimes c\mapsto h(f(c))].$
\\ 
$\gamma $ & $:$ & $A\otimes _{R}C^{\ast }$ & $\rightarrow $ & $\mathrm{Hom}%
_{R}(C,A),$ & $a\otimes g$ & $\mapsto $ & $[c\mapsto g(c)a].$%
\end{tabular}
\end{equation*}
By Corollary \ref{conseq} $(\mathrm{Hom}_{R}(C,A),A\star C)\in \mathcal{P}%
_{m}^{\alpha },$ $\gamma (A\otimes _{R}C^{\ast })\subset \mathrm{Hom}%
_{R}(C,A)$ is dense (w.r.t. the left $C$-adic topology) and we get category
isomorphisms 
\begin{equation*}
\begin{tabular}{lllll}
$\mathcal{M}^{A\star C}$ & $\simeq $ & $\mathrm{Rat}^{A\star C}(_{A\otimes
_{R}C^{\ast }}\mathcal{M})$ & $=$ & $\sigma \lbrack _{A\otimes _{R}C^{\ast
}}(A\star C)]$ \\ 
& $\simeq $ & $\mathrm{Rat}^{A\star C}(_{(A\star C)^{\ast }}\mathcal{M})$ & $%
=$ & $\sigma \lbrack _{(A\star C)^{\ast }}(A\star C)]$ \\ 
& $\simeq $ & $\mathrm{Rat}^{A\star C}(_{\mathrm{Hom}_{R}(C,A)}\mathcal{M})$
& $=$ & $\sigma \lbrack _{\mathrm{Hom}_{R}(C,A)}(A\star C)].$%
\end{tabular}
\end{equation*}
\end{punto}

\begin{proposition}
\label{A*C}If $R$ is noetherian, then we have bifunctors 
\begin{equation}
-\star -:\mathbf{Alg}_{R}^{\alpha }\times \mathbf{Cog}\rightarrow \mathbf{Cog%
}_{R}\text{ and }-\star -:\mathbf{Alg}_{R}^{\alpha }\times \mathbf{Cog}%
{\normalsize _{R}^{\alpha }}\rightarrow \mathbf{Cog}_{R}^{\alpha }.
\label{(-)^0x-}
\end{equation}
\end{proposition}

\begin{Beweis}
Let $A\in \mathbf{Alg}_{R}^{\alpha }.$ Then $A^{\circ }$ is by Theorem \ref
{R(G)-co} a locally projective $R$-coalgebra (i.e. an $\alpha $-coalgebra).
If $C$ is a (locally projective) $R$-coalgebra, then $A\star C:=A^{\circ
}\otimes _{R}C$ is a (locally projective) $R$-coalgebra by \cite[Lemma 2.8]
{Abu}). Analog to \cite{Par73} one can see that (\ref{(-)^0x-}) are
bifunctors.$\blacksquare $
\end{Beweis}

\subsection*{Continuous Dual Coalgebras}

\begin{definition}
Let $A$ be an $R$-algebra, $\mathcal{K}_{A}$ be the class of \emph{all }$R$%
-cofinite $A$-ideals and 
\begin{equation*}
\mathcal{E}_{A}:=\{I\vartriangleleft A\mid A/I\text{ is finitely generated
projective}\}.
\end{equation*}
For every subclass $\frak{F}\subseteq \mathcal{K}_{A}$ set 
\begin{equation*}
A_{\frak{F}}^{\circ }:=\{f\in A^{\ast }\mid f(I)=0\text{ for some }I\in 
\frak{F}\}.
\end{equation*}

\begin{enumerate}
\item  We call a filter $\frak{F}=\{I_{\lambda }\}_{\Lambda }$ consisting of 
$R$-cofinite $A$-ideals:

an $\alpha $\emph{-filter}, if the $R$-pairing $(A,A_{\frak{F}}^{\circ })$
satisfies the $\alpha $-condition;

\emph{cofinitary}, if $\frak{F}\cap \mathcal{E}_{A}$ is a filter basis of $%
\frak{F};$

\emph{cofinitely }$R$\emph{-cogenerated,} if $A/I$ is $R$-cogenerated for
every $I\in \frak{F}.$

\item  We call $A:$

an $\alpha $\emph{-algebra}, if $\mathcal{K}_{A}$ is an $\alpha $-filter;

\emph{cofinitary}, if $\mathcal{K}_{A}$ is a cofinitary filter;

\emph{cofinitely }$R$\emph{-cogenerated}, if $A/I$ is $R$-cogenerated for
every $I\in \mathcal{K}_{A}.$
\end{enumerate}
\end{definition}

\begin{definition}
(\cite{Tak81}) An $R$-coalgebra $C$ is called \emph{infinitesimal flat,}%
{\normalsize \ }if $C=\underrightarrow{lim}C_{\lambda }$ for a directed
system of \emph{finitely generated projective} $R$-subcoalgebras $%
\{C_{\lambda }\}_{\Lambda }.$
\end{definition}

\begin{proposition}
\label{A^0-F-co}Let $A$ be an $R$-algebra, $\frak{F}$ be a filter consisting
of $R$-cofinite $A$-ideals, $P:=(A,A_{\frak{F}}^{\circ })$ and consider $A$
as a left \emph{(}a right\emph{)} linear topological $R$-algebra with the
induced topology $\frak{T}(\frak{F}).$

\begin{enumerate}
\item  Assume $\frak{F}$ to be cofinitely $R$-cogenerated. Then $\frak{T}(%
\frak{F})$ is Hausdorff, iff $\kappa _{P}:A\rightarrow A_{\frak{\mathcal{F}}%
}^{\circ \ast }$ is an embedding.

\item  Assume $R$ to be noetherian and $\frak{F}$ to be an $\alpha $-filter.
Then $A_{\frak{F}}^{\circ }$ is an $\alpha $-coalgebra, $(A,A_{\frak{F}%
}^{\circ })\in \mathcal{P}_{m}^{\alpha }$ and $\kappa _{P}(A)\subset A_{%
\frak{F}}^{\circ \ast }$ is dense \emph{(}w.r.t. the finite topology\emph{)}.

\item  If $A/I$ is $R$-reflexive for every $I\in \frak{F}$ \emph{(}e.g. $R$
is an injective cogenerator\emph{)}, then $\widehat{A}\simeq A_{\frak{F}%
}^{\circ \ast }$ as left \emph{(}right\emph{)} linear topological $R$%
-modules.
\end{enumerate}
\end{proposition}

\begin{Beweis}
\begin{enumerate}
\item  By assumption $A/I$ is $R$-cogenerated for every $I\in \frak{F},$
hence 
\begin{equation*}
\overline{0_{A}}=\bigcap\limits_{I{\normalsize \in }\frak{F}}I=\bigcap_{I\in 
\frak{F}}\mathrm{KeAn}{\normalsize (I)}=\mathrm{Ke}(\sum\limits_{I\in \frak{F%
}}\mathrm{An}(I))=\mathrm{Ke}(A_{\frak{F}}^{\circ })=\mathrm{Ke}(\kappa
_{P}).
\end{equation*}

\item  Every $I\in \frak{F}\ $is a two-sided $A$-ideal and so $A_{\frak{F}%
}^{\circ }\subset A^{\circ }$ is an $A$-subbimodule. The results follows
then from Theorem \ref{R(G)-co}.

\item  If $A/I$ is $R$-reflexive for every $I\in \frak{F},$ then we have
isomorphisms of topological $R$-modules 
\begin{equation*}
\widehat{A}=\underleftarrow{lim}_{\frak{F}}A/I\simeq \underleftarrow{lim}_{%
\frak{F}}(A/I)^{\ast \ast }\simeq ({\underrightarrow{lim}}_{\frak{F}%
}(A/I)^{\ast })^{\ast }=:(A_{\frak{F}}^{\circ })^{\ast }.
\end{equation*}
If $R$ is an injective cogenerator, then all finitely generated $R$-modules
are $R$-reflexive (e.g. \cite[48.13]{Wis88}) and we are done.$\blacksquare $
\end{enumerate}
\end{Beweis}

The following result extends observations in \cite{Lar98} (resp. \cite
{AG-TL2001}) on \emph{cofinitary }$R$\emph{-algebras} over Dedekind domains
(resp. noetherian rings) to the case of \emph{cofinitary filters}\ for
algebras over arbitrary commutative base rings:

\begin{proposition}
\label{A-cof}Let $A$ be an $R$-algebra, $\frak{F}$ be a filter consisting of 
$R$-cofinite $A$-ideals, $P:=(A,A_{\frak{F}}^{\circ })$ and consider $A$ as
a left \emph{(}a right\emph{)} linear topological $R$-algebra with the
induced left \emph{(}right\emph{)} linear topology $\frak{T}(\frak{F}).$ If $%
\frak{F}\ $is cofinitary, then

\begin{enumerate}
\item  $\frak{T}(\frak{F})$ is Hausdorff, iff $\kappa _{P}:A\rightarrow A_{%
\frak{F}}^{\circ \ast }$ is an embedding.

\item  $A_{\frak{F}}^{\circ }$ is an infinitesimal flat $\alpha $-coalgebra, 
$P\in \mathcal{P}_{m}^{\alpha }$ and $\kappa _{P}(A)\subset A_{\frak{F}%
}^{\circ \ast }$ is dense.

\item  $\widehat{A}\simeq A_{\frak{F}}^{\circ \ast }$ as left \emph{(}right%
\emph{)} linear topological $R$-algebras.
\end{enumerate}
\end{proposition}

\begin{Beweis}
\begin{enumerate}
\item  For every $I\in \mathcal{E}_{A}$ the $R$-module $A/I\;$is in
particular $R$-cogenerated and the result follows from Proposition \ref
{A^0-F-co} (1).

\item  For $I,J\in \frak{F}\cap \mathcal{E}_{A}$ set $I\leq J,$ if $%
I\supseteq J$ and consider the canonical $R$-algebra epimorphism $\pi
_{I,J}:A/J\rightarrow A/I.$ Then 
\begin{equation*}
\{((A/I)^{\ast },\pi _{I,J}^{\ast })|\text{ }I\in \frak{F}\cap \mathcal{E}%
_{A},\text{ }\pi _{I,J}^{\ast }:(A/I)^{\ast }\hookrightarrow (A/J)^{\ast }\}
\end{equation*}
is a directed system of finitely generated projective $R$-coalgebras with $R$%
-coalgebra morphisms $\pi _{I,J}^{\ast }:(A/I)^{\ast }\rightarrow
(A/J)^{\ast }.$ Then $A_{\frak{F}}^{\circ }=A_{\frak{F}\cap \mathcal{E}%
_{A}}^{\circ }\simeq {\underrightarrow{lim}}_{\frak{F}\cap \mathcal{E}%
_{A}}(A/I)^{\ast }$ is an infinitesimal flat $R$-coalgebra.

Let $M$ be an arbitrary $R$-module. If $\sum\limits_{i=1}^{k}m_{i}\otimes
g_{i}\in \mathrm{Ke}(\alpha _{M}^{P}),$ then there exists $I\in \frak{F}\cap 
\mathcal{E}_{A},$ such that $\{g_{1},...,g_{n}\}\subset \mathrm{An}(I).$ If $%
\{(a_{l}+I,f_{l})\}_{l=1}^{k}$ is a dual basis for $(A/I)^{\ast },$ then 
\begin{equation*}
\begin{tabular}{lll}
$\sum\limits_{i=1}^{n}m_{i}\otimes g_{i}$ & $=$ & $\sum%
\limits_{i=1}^{n}m_{i}\otimes (\sum\limits_{l=1}^{k}g_{i}(a_{l}+I)f_{l})$ \\ 
& $=$ & $\sum\limits_{i=1}^{n}m_{i}\otimes
(\sum\limits_{l=1}^{k}g_{i}(a_{l})f_{l})$ \\ 
& $=$ & $\sum\limits_{l=1}^{k}(\sum\limits_{i=1}^{n}g_{i}(a_{l})m_{i})%
\otimes f_{l}=0.$%
\end{tabular}
\end{equation*}
Obviously the canonical $R$-linear mapping $\kappa _{P}:A\rightarrow A_{%
\frak{F}}^{\circ \ast }$ is an $R$-algebra morphism, i.e. $P$ is a measuring 
$\alpha $-pairing. The density of $\kappa _{P}(A)\subset A_{\frak{F}}^{\circ
\ast }$ follows then by Theorem \ref{R(G)-co}.

\item  For every $I\in \frak{F}\cap \mathcal{E}_{A}$ the $R$-module $A/I$ is
finitely generated projective, hence $(A/I)^{\ast }$ is an $R$-coalgebra and 
$(A/I)^{\ast \ast }\simeq A/I$ as $R$-algebras. So we have an isomorphisms
of topological $R$-algebras 
\begin{equation*}
\widehat{A}=\underleftarrow{lim}_{\frak{F}\cap \mathcal{E}_{A}}A/I\simeq 
\underleftarrow{lim}_{\frak{F}\cap \mathcal{E}_{A}}(A/I)^{\ast \ast }\simeq (%
{\underrightarrow{lim}}_{\frak{F}\cap \mathcal{E}_{A}}(A/I)^{\ast })^{\ast
}=({\underrightarrow{lim}}_{\frak{F}}(A/I)^{\ast })^{\ast }=(A_{\frak{F}%
}^{\circ })^{\ast }.\blacksquare
\end{equation*}
\end{enumerate}
\end{Beweis}

\qquad As a consequence of Propositions \ref{A^0-F-co}, \ref{A-cof} and
Theorem \ref{cor-dicht} we get

\begin{corollary}
Let $A$ be an $R$-algebra and $\frak{F}\ $be a filter consisting of $R$%
-cofinite $A$-ideals. If $R$ is noetherian and $\frak{F}$ is an $\alpha $%
-filter, or $\frak{F}$ is cofinitary, then we have isomorphisms of
categories 
\begin{equation*}
\begin{tabular}{lllll}
$\mathcal{M}^{A_{\frak{F}}^{\circ }}$ & $\simeq $ & $\mathrm{Rat}^{A_{\frak{F%
}}^{\circ }}(_{A}\mathcal{M})$ & $=$ & $\sigma \lbrack _{A}A_{\frak{F}%
}^{\circ }]$ \\ 
& $\simeq $ & $\mathrm{Rat}^{A_{\frak{F}}^{\circ }}(_{A_{\frak{F}}^{\circ
\ast }}\mathcal{M})$ & $=$ & $\sigma \lbrack _{A_{\frak{F}}^{\circ \ast }}A_{%
\frak{F}}^{\circ }]$%
\end{tabular}
\text{ }\&\text{ } 
\begin{tabular}{lllll}
$^{A_{\frak{F}}^{\circ }}\mathcal{M}$ & $\simeq $ & $^{A_{\frak{F}}^{\circ }}%
\mathrm{Rat}(\mathcal{M}_{A})$ & $=$ & $\sigma \lbrack A_{\frak{F}}^{\circ
}{}_{A}]$ \\ 
& $\simeq $ & $^{A_{\frak{F}}^{\circ }}\mathrm{Rat}(\mathcal{M}_{A_{\frak{F}%
}^{\circ \ast }})$ & $=$ & $\sigma \lbrack A_{\frak{F}}^{\circ }{}_{A_{\frak{%
F}}^{\circ \ast }}].$%
\end{tabular}
\end{equation*}
\end{corollary}

\begin{punto}
\label{filt-2}Let $A,B$ be $R$-algebras, $\frak{F}_{A},$ $\frak{F}_{B}$ be
filter bases consisting of $R$-cofinite $A$-ideals, $B$-ideals respectively
and 
\begin{equation}
\frak{F}_{A}\times \frak{F}_{B}:=\{\mathrm{Im}(\iota _{I}\otimes id_{B})+%
\mathrm{Im}(id_{A}\otimes \iota _{J})|\text{ }I\in \frak{F}_{A},\text{ }J\in 
\frak{F}_{B}\}.  \label{filt-mal}
\end{equation}
Obviously $\frak{F}_{A}\times \frak{F}_{B}$ is a filter basis consisting of $%
R$-cofinite $A\otimes _{R}B$-ideals and induces so a linear topology $\frak{T%
}(\frak{F}_{A}\times \frak{F}_{B})$ on $A\otimes _{R}B,$ such that $%
(A\otimes _{R}B,\frak{T}(\frak{F}_{A}\times \frak{F}_{B}))$ is a linear
topological $R$-algebra and $\frak{F}_{A}\times \frak{F}_{B}\ $is a
neighbourhood basis of $0_{A\otimes _{R}B}.$
\end{punto}

\qquad One can generalize \cite[Proposition 4.9, Theorem 4.10]{AG-TL2001} to
the following more general result:

\begin{theorem}
\label{pro-2}Let $A,B$ be $R$-algebras, $\frak{F}_{A},$ $\frak{F}_{B}$ be
filters consisting of $R$-cofinite $A$-deals, $B$-ideals respectively and
consider the canonical $R$-linear mapping $\delta :A^{\ast }\otimes
_{R}B^{\ast }\rightarrow (A\otimes _{R}B)^{\ast }.$

\begin{enumerate}
\item  If $\frak{F}_{A}\ $and $\frak{F}_{B}$ are cofinitary, then the filter
of $A\otimes _{R}B$-ideals generated by $\frak{F}_{A}\times \frak{F}_{B}$ is
cofinitary and $(A\otimes _{R}B)_{\frak{F}_{A}\times \frak{F}_{B}}^{\circ }$
is an $R$-coalgebra. If $R$ is noetherian, then $\delta $ induces an $R$%
-coalgebra isomorphism 
\begin{equation*}
A_{\frak{F}{\normalsize _{A}}}^{\circ }\otimes _{R}B_{\frak{F}{\normalsize %
_{B}}}^{\circ }\simeq (A\otimes _{R}B)_{\frak{F}{\normalsize _{A}\times }%
\frak{F}{\normalsize _{B}}}^{\circ }.
\end{equation*}

\item  Let $R\;$be noetherian. If $\frak{F}_{A}$ is an $\alpha $-filter and $%
\frak{F}_{B}$ is cofinitary, then the filter generated by $\frak{F}%
_{A}\times \frak{F}_{B}$ is an $\alpha $-filter, $(A\otimes _{R}B)_{\frak{F}%
_{A}\times \frak{F}_{B}}^{\circ }$ is an $R$-coalgebra and $\delta $ induces
an $R$-coalgebra isomorphism 
\begin{equation*}
A_{\frak{F}{\normalsize _{A}}}^{\circ }\otimes _{R}B_{\frak{F}{\normalsize %
_{B}}}^{\circ }\simeq (A\otimes _{R}B)_{\frak{F}{\normalsize _{A}\times }%
\frak{F}{\normalsize _{B}}}^{\circ }.
\end{equation*}
\end{enumerate}
\end{theorem}

\begin{theorem}
\label{erb}Let $R$ be hereditary and noetherian.

\begin{enumerate}
\item  All $R$-algebras satisfy the $\alpha $-condition, i.e. $\mathbf{Alg}%
_{R}^{\alpha }=\mathbf{Alg}_{R}.$

\item  There is a \emph{duality }between{\normalsize \ }$\mathbf{Alg}_{R}$
and $\mathbf{Cog}_{R}$ through the right-adjoint contravariant functors 
\begin{equation*}
(-)^{\ast }:\mathbf{Cog}_{R}\rightarrow \mathbf{Alg}_{R},\text{ }(-)^{\circ
}:\mathbf{Alg}_{R}\rightarrow \mathbf{Cog}_{R}.
\end{equation*}
\end{enumerate}
\end{theorem}

\begin{Beweis}
\begin{enumerate}
\item  Let $A$ be an arbitrary $R$-algebra. By \cite[Proposition 2.11]
{AG-TW2000} $A^{\circ }\subset R^{A}$ is pure and so $(A,A^{\circ })$ is an $%
\alpha $-pairing by \cite[Proposition 2.6]{Abu}.

\item  For every $R$-algebra $A$ the canonical mapping $\lambda
_{A}:A\rightarrow A^{\circ \ast }$ is an $R$-algebra morphism and for every $%
R$-coalgebra $C$ the canonical mapping $\Phi _{C}:C\rightarrow C^{\ast \circ
}$ is an $R$-coalgebra morphism (compare Lemma \ref{alg-coalg}). Moreover
for every $A\in \mathbf{Alg}_{R}$ and every $C\in \mathbf{Cog}_{R}$%
\begin{equation*}
\Upsilon _{A,C}:\mathrm{Alg}_{R}(A,C^{\ast })\rightarrow \mathrm{Cog}%
_{R}(C,A^{\circ }),\text{ }\xi \mapsto \xi ^{\circ }\circ \Phi _{C}
\end{equation*}
is an isomorphism with inverse 
\begin{equation*}
\Psi _{A,C}:\mathrm{Cog}_{R}(C,A^{\circ })\rightarrow \mathrm{Alg}%
_{R}(A,C^{\ast }),\text{ }\theta \mapsto \theta ^{\ast }\circ \lambda _{A}.
\end{equation*}
It's easy to see that $\Upsilon _{A,C}$ and $\Psi _{A,C}$ are functorial in $%
A$ and $C.\blacksquare $
\end{enumerate}
\end{Beweis}

\subsection*{Locally finite modules}

\begin{punto}
\textbf{The cofinite topology.}\label{ko-top}{\normalsize \ }Let $R$ be
noetherian. For every $R$-algebra $A,$ the class $\mathcal{K}_{A}$ of \emph{%
all }$R$-cofinite \emph{two-sided }$A$-ideals is obviously a filter and
induces on $A$ a left (a right) linear topology, the so called \emph{left\ }(%
\emph{right}) \emph{cofinite topology }$\mathrm{Cf}(A),$\textbf{\ }such that 
$\mathcal{K}_{A}$ is neighbourhood basis of $0_{A}.$ If $A^{\circ }:=A_{%
\mathcal{K}_{A}}^{\circ }$ is an $R$-coalgebra, then we call it the \emph{%
continuous dual }$R$\emph{-coalgebra} of $A.$

Consider $A$ with the left linear topology $\mathrm{Cf}(A).$ Let $M$ be a
left $A$-module and consider the filter $\mathcal{K}_{M}$ of \emph{all }$R$%
-cofinite $A$-submodules of $M.$ Let $L\subset M$ be an $R$-cofinite $A$%
-submodule and consider the $R$-linear mapping 
\begin{equation*}
\varphi _{L}:A\rightarrow \mathrm{End}_{R}(M/L),\text{ }a\mapsto \lbrack
m+L\mapsto am+L].
\end{equation*}
Then $A/\mathrm{Ke}(\varphi _{L})\hookrightarrow \mathrm{End}_{R}(M/L)$ and
so 
\begin{equation*}
I_{L}:=\mathrm{Ke}(\varphi _{L})=\{a\in A|\text{ }aM\subseteq L\}
\end{equation*}
is an $R$-cofinite two-sided $A$-ideal. If $m\in M$ is arbitrary, then $%
I_{L}:=(L:M)\subset (L:m),$ hence $(L:m)$ is open w.r.t. to the left
cofinite topology $\mathrm{Cf}(A).$ So $M$ becomes a topology, the so called 
\emph{cofinite topology} $\mathrm{Cf}(M),$ such that $(M,\mathrm{Cf}(M))$ is
a linear topological left $(A,\mathrm{Cf}(A))$-module and $\mathcal{K}_{M}$
is a neighbourhood basis of $0_{M}.$

Considering $A$ with the right cofinite topology $\mathrm{Cf}(A)$ it turns
out that for every right $A$-module $M,$ the filter of $R$-cofinite right $A$%
-submodules of $M$ induces on $M$ a topology, the cofinite topology $\mathrm{%
Cf}(M),$ such that $(M,\mathrm{Cf}(M))$ is a linearly topological right $(A,%
\mathrm{Cf}(A))$-module.
\end{punto}

\begin{punto}
Let $R$ be noetherian and $A$ be an $R$-algebra. A left $A$-module $M$ is
called \emph{locally finite}, if $Am$ is finitely generated for every $m\in
M.$ With $\mathrm{Loc}(_{A}\mathcal{M})\subset $ $_{A}\mathcal{M}$ we denote
the \emph{full }subcategory of locally finite left $A$-modules. For every
left $A$-module $M$ it follows that $\mathrm{Loc}(M)\subset M$ is an $A$%
-submodule (since the ground ring $R$ is noetherian) and we get a preradical 
\begin{equation*}
\mathrm{Loc}(-):\text{ }_{A}\mathcal{M}\rightarrow {\normalsize _{A}}%
\mathcal{M},\text{ }M\mapsto \{m\in M|\text{ }Am\text{ is finitely generated
in }\mathcal{M}_{R}\}
\end{equation*}
with pretorsion class the full subcategory $\mathrm{Loc}(_{A}\mathcal{M}%
)\subseteq $ $_{A}\mathcal{M}$ of locally finite left $A$-modules.

Analogously one defines the preradical $\mathrm{Loc}(-):\mathcal{M}%
_{A}\rightarrow \mathcal{M}_{A}$ with pretorsion class the full subcategory
of locally finite right $A$-modules $\mathrm{Loc}(\mathcal{M}_{A})\subseteq 
\mathcal{M}_{A}.$
\end{punto}

\begin{lemma}
\label{grun}Let $R\;$be noetherian and $A$ be an $R$-algebra. For every
right $A$-module $M$ we have 
\begin{equation}
\begin{tabular}{lll}
$M^{\circ }$ & $:=$ & $\{f\in M^{\ast }|$ $f(MI)=0$ for some $R$-cofinite 
\emph{(}right\emph{)} $A$-ideal $I\subset A\}$ \\ 
& $=$ & $\{f\in M^{\ast }|$ $Af$ is finitely generated in $\mathcal{M}_{R}\}$
\emph{(}$=\mathrm{Loc}(_{A}M^{\ast })$\emph{)} \\ 
& $=$ & $\{f\in M^{\ast }|\text{ }f(L)=0\text{ for some}$\ $R$-cofinite
right $A$-submodule $L\subset M\}.$%
\end{tabular}
\label{M^0=}
\end{equation}
\end{lemma}

\begin{Beweis}
Let $f\in M^{\ast }$ with $f(MI)=0$ for an $R$-cofinite right $A$-ideal $I.$
If $\{a_{1}+I,...,a_{k}+I\}$ is a generating system for $A/I$ over $R,$ then 
$\{a_{1}f,...,a_{k}f\}$ is a generating system for $Af$ over $R,$ i.e. $f\in 
\mathrm{Loc}(_{A}M^{\ast }).$

Let $f\in \mathrm{Loc}(_{A}M^{\ast })$ and assume that $Af=\sum%
\limits_{i=1}^{k}Rf_{i}$ with $\{f_{1},...,f_{k}\}\subset M^{\ast }.$ Then $%
L:=\mathrm{Ke}(Af)=\bigcap\limits_{i=1}^{k}\mathrm{Ke}(f_{i})\subset M$ is a
right $A$-submodule and moreover $M/L\hookrightarrow
\bigoplus\limits_{i=1}^{k}M/\mathrm{Ke}(f_{i}),$ i.e. $L\subset M$ is an $R$%
-cofinite $A$-submodule.

Let $f\in (M/L)^{\ast }\simeq \mathrm{An}(L)$ for some $R$-cofinite $A$%
-submodule $L\subseteq M.$ Then $I_{L}:=(L:M)\;$is an $R$-cofinite \emph{%
two-sided} $A$-ideal (compare \ref{ko-top}) and moreover $f(MI_{L})\subseteq
f(L)=0,$ i.e. $f\in M^{\circ }.\blacksquare $
\end{Beweis}

\qquad It's well known that for an $R$-algebra $A$ over a base field $R,$
the category of right (left) $A^{\circ }$-comodules and the category of
locally finite left (right) $A$-modules coincide (e.g. \cite{Abe80}, \cite
{Rad73}). Over arbitrary commutative rings we have

\begin{proposition}
\label{A^0-cor}Let $R$ be noetherian and $A$ be an $R$-algebra.

\begin{enumerate}
\item  Every $A^{\circ }$-subgenerated left (right) $A$-module is locally
finite.

\item  If $A$ is cofinitely $R$-cogenerated, then $\sigma \lbrack
_{A}A^{\circ }]=\mathrm{Loc}(_{A}\mathcal{M})$ and $\sigma \lbrack
A_{A}^{\circ }]=\mathrm{Loc}(\mathcal{M}_{A}).$

\item  If $A$ is an \emph{(}cofinitely $R$-cogenerated\emph{)} $\alpha $%
-algebra, then we have category isomorphisms 
\begin{equation*}
\begin{tabular}{llll}
$\mathcal{M}^{A^{\circ }}$ & $\simeq $ & $\mathrm{Rat}^{A^{\circ }}(_{A}%
\mathcal{M})$ & $=\sigma \lbrack _{A}A^{\circ }]$ \\ 
& $\simeq $ & $\mathrm{Rat}^{A^{\circ }}(_{A^{\circ \ast }}\mathcal{M})$ & $%
=\sigma \lbrack _{A^{\circ \ast }}A^{\circ }]$%
\end{tabular}
\&
\begin{tabular}{llll}
$^{A^{\circ }}\mathcal{M}$ & $\simeq $ & $^{A^{\circ }}\mathrm{Rat}(\mathcal{%
M}_{A})$ & $=\sigma \lbrack A_{A}^{\circ }]$ \\ 
& $\simeq $ & $^{A^{\circ }}\mathrm{Rat}(\mathcal{M}_{A^{\circ \ast }})$ & $%
=\sigma \lbrack A_{A^{\circ \ast }}^{\circ }].$%
\end{tabular}
\end{equation*}
If $A$ is moreover cofinitely $R$-cogenerated, then 
\begin{equation*}
\mathcal{M}^{A^{\circ }}\simeq \mathrm{Loc}(_{A}\mathcal{M})\text{ and }%
^{A^{\circ }}\mathcal{M}\simeq \mathrm{Loc}(\mathcal{M}_{A}).
\end{equation*}
\end{enumerate}
\end{proposition}

\begin{Beweis}
\begin{enumerate}
\item  Let $M\in \sigma \lbrack _{A}A^{\circ }].$ Then there exists for
every $m\in M$ a finite subset $W=\{f_{1},...,f_{k}\}\subset A^{\circ },$
such that $(0_{A^{\circ }}:W)\subset (0_{M}:m).$ Choose for every $i=1,...,k$
an $R$-cofinite $A$-ideal $J_{i}\subset Ke(f_{i})$ and consider the $R$%
-cofinite $A$-ideal $J:=\bigcap_{i=1}^{k}J_{i}.$ If $a\in J,$ then for every 
$\widetilde{a}\in A$ and $i=1,...,k$ we have $(a\rightharpoonup f_{i})(%
\widetilde{a})=f_{i}(\widetilde{a}a)=0.$ Consequently $J\subset (0_{A^{\circ
}}:W)\subset (0_{M}:m)$ and so $Am\simeq A/(0_{M}:m)\;$is finitely generated
in $\mathcal{M}_{R}.$ Hence $_{A}M$ is locally finite.

\item  By (1)\ $\sigma \lbrack _{A}A^{\circ }]\subset \mathrm{Loc}(_{A}%
\mathcal{M}).$ Assume now that $A$ is cofinitely $R$-cogenerated. Let $N$ be
a locally finite left $A$-module. For every $n\in N$ the $R$-module $%
A/(0_{N}:n)\simeq An$ is finitely generated and so there exists by Remark 
\ref{twoside} an $R$-cofinite $A$-ideal $I\subset (0_{N}:n).$ By assumption $%
A/I$ is $R$-cogenerated and so $I=\mathrm{KeAn}(I)$ (e.g. \cite[28.1]{Wis88}%
). If $\mathrm{An}(I)\simeq (A/I)^{\ast }=\sum\limits_{i=1}^{k}Rg_{i}$ and $%
W:=\{g_{1},...,g_{k}\},$ then it follows for every $a\in (0_{A^{\circ }}:W)$
that $g_{i}(a)=(a\rightharpoonup g_{i})(1_{A})=0.$ So $(0_{A^{\circ
}}:W)\subset \mathrm{KeAn}(I)=I\subset (0_{N}:n),$ i.e. $_{A}N$ is $A^{\circ
}$-subgenerated.

\item  The category isomorphisms follow from Theorem \ref{cor-dicht} (and
(2)).$\blacksquare $
\end{enumerate}
\end{Beweis}

\section{Dual comodules}

\qquad

\qquad In this section we discuss for every $(A,C)\in \mathcal{P}%
_{m}^{\alpha }$ the \emph{duality} between the category of right (left) $A$%
-modules and the category the right (left) $C$-comodules.

\begin{punto}
\label{box}Let $P=(A,C)\in \mathcal{P}_{m}.$ By Theorem \ref{cor-dicht} $%
\mathcal{M}^{C}\subset $ $_{A}\mathcal{M}$ is a subcategory and so we have a
contravariant functor 
\begin{equation}
(-)^{\ast }:\mathcal{M}^{C}\rightarrow \mathcal{M}_{A},\text{ }(N,\varrho
_{N})\mapsto (N^{\ast },\rho _{N^{\ast }}),  \label{(-)^*}
\end{equation}
where 
\begin{equation}
\rho _{N^{\ast }}:N^{\ast }\rightarrow \mathrm{Hom}_{R}(A,N^{\ast }),\text{ }%
f\mapsto \lbrack a\mapsto \lbrack n\mapsto \sum f(n_{<0>})<a,n_{<1>}>]].
\label{vp_nst}
\end{equation}
If moreover $P$ satisfies the $\alpha $-condition, then we get by Theorem 
\ref{cor-dicht} the contravariant functor 
\begin{equation*}
(-)^{r}:\mathcal{M}_{A}\rightarrow \mathcal{M}^{C},\text{ }M\mapsto M^{r}:=%
\mathrm{Rat}^{C}(_{A}M^{\ast }).
\end{equation*}
If $M,$ $\widetilde{M}$ are right $A$-modules and $f\in \mathrm{Hom}_{A-}(M,%
\widetilde{M}),$ then $f^{\ast }\in \mathrm{Hom}_{-A}(N^{\ast },\widetilde{M}%
^{\ast })$ and we denote with $f^{r}\in \mathrm{Hom}^{C}(\widetilde{M}%
^{r},M^{r})$ the restriction of $f^{\ast }$ to $\widetilde{M}^{r}\subset 
\widetilde{M}^{\ast }$ (see Lemma \ref{clos} (4) and Theorem \ref{cor-dicht}%
). For every right $A$-module $M$ we call $M^{r}$ the \emph{dual }$C$\emph{%
-comodule }of $M$\ w.r.t. $P.$
\end{punto}

\begin{punto}
\label{fp}(\cite{Wis2000}) Let $_{R}C\;$be a flat $R$-coalgebra, $N$ be a
left $C$-comodule and consider the $R$-linear mapping 
\begin{equation}
\gamma :N^{\ast }\rightarrow \mathrm{Hom}_{R}(N,C),\text{ }f\mapsto \lbrack
n\mapsto \sum f(n_{<-1>})n_{<0>}].  \label{rh_Mst}
\end{equation}
If $_{R}N$ is finitely presented, then by $N^{\ast }\otimes _{R}C\simeq 
\mathrm{Hom}_{R}(N,C)$ (e.g. \cite[15.7]{Wis96}) and $N^{\ast }$ becomes a
structure of a right $C$-comodule through 
\begin{equation}
\varrho _{M^{\ast }}:N^{\ast }\overset{\gamma }{\rightarrow }\mathrm{Hom}%
_{R}(N,C)\simeq N^{\ast }\otimes _{R}C.  \label{vr_Mst}
\end{equation}
Analog, if $N\in \mathcal{M}^{C}$ and $N_{R}$ is finitely presented, then $%
N^{\ast }$ becomes a left $C$-comodule structure.
\end{punto}

\begin{theorem}
\label{adj}For every $(A,C)\in \mathcal{P}_{m}^{\alpha }\;$there is a \emph{%
duality} between the category of right $C$-comodules and the category of
right $A$-modules through the right adjoint contravariant functors 
\begin{equation*}
(-)^{\ast }:\mathcal{M}^{C}\rightarrow \mathcal{M}_{A}\text{ and }(-)^{r}:%
\mathcal{M}_{A}\rightarrow \mathcal{M}^{C}.
\end{equation*}
\end{theorem}

\begin{Beweis}
For every right $C$-comodule $N\;$the canonical mapping $\Phi
_{N}:N\rightarrow N^{\ast \ast }$ is $A$-linear, we have by Lemma \ref{clos}
(4) $\Phi _{N}(N)\subseteq N^{\ast r}$ and so $\Phi _{N}:N\rightarrow
N^{\ast r}$ is $C$-colinear by Theorem \ref{cor-dicht}. On the other hand,
for every right $A$-module $M$ the canonical mapping $\lambda
_{M}:M\rightarrow M^{r\ast }$ is $A$-linear. It's easy to see then that we
have functorial homomorphisms (in $M\in \mathcal{M}_{A}$ and $N\in \mathcal{M%
}^{C}$) 
\begin{equation*}
\begin{tabular}{llllllll}
$\Upsilon _{N,M}$ & $:$ & $\mathrm{Hom}_{-A}(M,N^{\ast })$ & $\rightarrow $
& $\mathrm{Hom}^{C}(N,M^{r}),$ & $f$ & $\mapsto $ & $f^{r}\circ \Phi _{N},$
\\ 
$\Psi _{N,M}$ & $:$ & $\mathrm{Hom}^{C}(N,M^{r})$ & $\rightarrow $ & $%
\mathrm{Hom}_{-A}(M,N^{\ast }),$ & $g$ & $\mapsto $ & $g^{\ast }\circ
\lambda _{M}.$%
\end{tabular}
\end{equation*}
Moreover $\Upsilon _{N,M}$ is bijective with inverse $\Psi
_{N,M}.\blacksquare $
\end{Beweis}

\begin{Notation}
For every $R$-algebra $A$ denote with $\mathcal{M}_{A}^{f}$ ($_{A}\mathcal{M}%
^{f}$) the category of \emph{finitely generated}{\normalsize \ }right (left) 
$A$-modules.
\end{Notation}

\begin{lemma}
\label{dual-f}Let $R$ be noetherian. For every $(A,C)\in \mathcal{P}%
_{m}^{\alpha }$ there is a \emph{duality}\ between $\mathrm{Rat}^{C}(_{A}%
\mathcal{M}^{f})$ and $^{C}\mathrm{Rat}(\mathcal{M}_{A}^{f})$ through the
right-adjoint contravariant functors 
\begin{equation*}
(-)^{\ast }:\text{ }^{C}\mathrm{Rat}(\mathcal{M}_{A}^{f})\rightarrow \mathrm{%
Rat}^{C}(_{A}\mathcal{M}^{f})\text{ and }(-)^{\ast }:\mathrm{Rat}^{C}(_{A}%
\mathcal{M}^{f})\rightarrow \text{ }^{C}\mathrm{Rat}(\mathcal{M}_{A}^{f})%
\text{ }.
\end{equation*}
\end{lemma}

\begin{Beweis}
Let $M\in $ $^{C}\mathrm{Rat}(\mathcal{M}_{A}^{f})$ ($M\in $ $\mathrm{Rat}%
^{C}(_{A}\mathcal{M}^{f})$). By \cite[Folgerung 2.2.24]{Abu2001} every
finitely generated $C$-rational left $A$-module is finitely generated over $%
R,$ hence $M_{R}$ is finitely generated and so $_{A}M^{\ast }$ ($M_{A}^{\ast
}$) is finitely generated. By assumption $R$ is noetherian and so $M_{R}$ is
finitely presented. Consequently $M^{\ast }$ is by \ref{fp} a $C$-rational
left (right) $A$-module. The claimed duality follows then from Theorem \ref
{adj}.$\blacksquare $
\end{Beweis}

\begin{punto}
\label{C-C*-adj}If $C$ is a locally projective $R$-coalgebra, then we get by
Theorem \ref{adj} right-adjoint contravariant functors 
\begin{equation*}
\begin{tabular}{llllllll}
$(-)^{\ast }$ & $:$ & $\mathcal{M}^{C}$ & $\rightarrow $ & $\mathcal{M}%
_{C^{\ast }},$ & $N$ & $\mapsto $ & $N^{\ast },$ \\ 
$(-)^{\Box }$ & $:$ & $\mathcal{M}_{C^{\ast }}$ & $\rightarrow $ & $\mathcal{%
M}^{C},$ & $M$ & $\mapsto $ & $M^{\Box }:=\mathrm{Rat}^{C}(_{C^{\ast
}}M^{\ast }).$%
\end{tabular}
\end{equation*}
\end{punto}

\begin{lemma}
Let $R$ be an injective cogenerator and $C$ be a locally projective $R$%
-coalgebra. If $M$ a right $C^{\ast }$-module, $L\subset M$ is a $C^{\ast }$%
-submodule and $M^{\Box }\subset M^{\ast }$ is dense, then $L^{\Box }\subset
L^{\ast }$ is dense.
\end{lemma}

\begin{Beweis}
By Lemma \ref{clos} (4) $\iota _{L}^{\ast }(M^{\Box })\subset L^{\Box }$ and
it follows from \cite[Proposition 1.10 (3.b)]{Abu} that $\overline{\iota
_{L}^{\ast }(M^{\Box })}=\iota _{L}^{\ast }(\overline{M^{\Box }})=\iota
_{L}^{\ast }(M^{\ast })=L^{\ast }.\blacksquare $
\end{Beweis}

M. Takeuchi \cite{Tak74} studied the category of locally finite modules of a
commutative algebra over a base field. In what follows we transfer some
results obtained by him to the category \textrm{Rat}$^{C}(_{A}\mathcal{M})$
corresponding to a measuring $\alpha $-pairing $P=(A,C)\in \mathcal{P}%
_{m}^{\alpha }$ with $A$ a \emph{commutative }algebra over an arbitrary
commutative ground ring.

\begin{proposition}
\label{hom=box}Let $P=(A,C)\in \mathcal{P}_{m}^{\alpha }$ with $A$ \emph{%
commutative} and denote with $\mathcal{M}_{A}^{f}\subset $ $_{A}\mathcal{M}$
the full subcategory of \emph{finitely generated} $A$-modules. Then we have
an isomorphism of functors 
\begin{equation*}
{\normalsize \mathrm{Hom}_{A}(-,C)\simeq (-)^{r}}:\mathcal{M}{\normalsize %
_{A}^{f}}\rightarrow \mathcal{M}^{C}.
\end{equation*}
\end{proposition}

\begin{Beweis}
\textbf{Step 1.}{\normalsize \ }Let $M\in \mathcal{M}_{A}^{f}$ be arbitrary
and consider $\mathrm{Hom}_{A}(M,C)$ with the canonical $A$-module structure
induced by $M_{A}.$ For arbitrary $f\in \mathrm{Hom}_{A}(M,C)$ the $A$%
-subbimodule $N:=f(M)\subset C$ is by Theorem \ref{cor-dicht} a $C$%
-bicoideal. Moreover $N_{A}$ is finitely generated and so finitely generated
in $\mathcal{M}_{R}$ (since $\mathrm{Rat}^{C}(_{A}\mathcal{M})$ is by 
\cite[Folgerung 2.2.24]{Abu2001} $(A,R)$-finite). Assume that $%
N=\sum\limits_{i=1}^{l}Rc_{i}$ with $\Delta
(c_{i})=\sum\limits_{j=1}^{l_{i}}c_{ij}\otimes \widetilde{c}_{ij}$ for every 
$i=1,...,k$ and set $K:=\sum\limits_{i=1}^{l}\sum%
\limits_{j=1}^{l_{i}}Rc_{ij}.$ Then $K^{\bot }\subset (0:f)$ and so $f$ is
by Proposition \ref{n_rat} $C$-rational. By our choice $f\in \mathrm{Hom}%
_{A}(M,C)$ is arbitrary, i.e. $\mathrm{Hom}_{A}(M,C)\in \mathrm{Rat}^{C}(_{A}%
\mathcal{M}).$

\textbf{Step 2. }$(-)^{r}\simeq \mathrm{Hom}_{A}(-,C).$

Let $N\in \mathcal{M}^{C},$ $M\in \mathcal{M}_{A}$ and consider the $C$%
-comodule $\mathrm{Hom}_{A}(M,C).$ The result follows then from the
functorial isomorphisms: 
\begin{equation*}
\begin{tabular}{llll}
$\mathrm{Hom}^{C}(N,\mathrm{Hom}_{A}(M,C))$ & $=$ & $\mathrm{Hom}_{A}(N,%
\mathrm{Hom}_{A}(M,C))$ & (Theorem \ref{cor-dicht}) \\ 
& $\simeq $ & $\mathrm{Hom}_{A}(M,\mathrm{Hom}_{A}(N,C))$ &  \\ 
& $\simeq $ & $\mathrm{Hom}_{A}(M,\mathrm{Hom}^{C}(N,C))$ & (Theorem \ref
{cor-dicht}) \\ 
& $\simeq $ & $\mathrm{Hom}_{A}(M,N^{\ast })$ & (\ref{M*=End}) \\ 
& $\simeq $ & $\mathrm{Hom}^{C}(N,M^{r})$ & (Theorem \ref{adj}).$%
\blacksquare $%
\end{tabular}
\end{equation*}
\end{Beweis}

\qquad As a consequence of Proposition \ref{hom=box} we get

\begin{corollary}
\label{A0-inj}

\begin{enumerate}
\item  Let $R$ be noetherian, $A$ be an $\alpha $-algebra and consider the
functor 
\begin{equation*}
(-)^{0}=\mathrm{Rat}^{A^{\circ }}(-)\circ (-)^{\ast }:\mathcal{M}%
_{A}\rightarrow \mathcal{M}^{A^{\circ }},\text{ }M\mapsto M^{0}:=\mathrm{Rat}%
^{A^{\circ }}(_{A}M^{\ast })
\end{equation*}
If $A$ is \emph{commutative}, then we have a functorial isomorphism 
\begin{equation*}
\mathrm{Hom}_{A}(-,A^{\circ })\simeq (-)^{0}:\mathcal{M}_{A}^{f}\rightarrow 
\mathcal{M}^{A^{\circ }}.
\end{equation*}

\item  If $C$ is a \emph{cocommutative}{\normalsize \ }locally projective $R$%
-coalgebra, then we have a functorial isomorphism 
\begin{equation*}
\mathrm{Hom}_{C^{\ast }}(-,C)\simeq (-)^{\Box }:\mathcal{M}_{C^{\ast
}}^{f}\rightarrow \mathcal{M}^{C}.
\end{equation*}
\end{enumerate}
\end{corollary}

\begin{corollary}
\label{^r-exact}Let $P=(A,C)\in \mathcal{P}_{m}^{\alpha },$ where $A$ is
commutative and noetherian. If $(-)^{r}:\mathcal{M}_{A}^{f}\rightarrow
\sigma \lbrack _{A}C]$ is exact, then $C$ is a injective $A$-module.
\end{corollary}

\begin{Beweis}
By \emph{Baer's criteria} (e.g. \cite[16.4]{Wis88}) it's enough to show that 
$C$ is $A$-injective. Let $I$ be an $A$-ideal. Then $I_{A}$ is finitely
generated and by assumption the following set mapping is surjective 
\begin{equation*}
A^{r}\overset{\iota ^{r}}{\rightarrow }I^{r}\rightarrow 0.
\end{equation*}
By Proposition \ref{hom=box} $\mathrm{Hom}_{A}(-,C)\simeq (-)^{r}$ and so 
\begin{equation*}
\mathrm{Hom}_{A}(A,C)\overset{(\iota ,C)}{\rightarrow }\mathrm{Hom}%
_{A}(I,C)\rightarrow 0
\end{equation*}
is a surjective set mapping, i.e. $C$ is $A$-injective and we are done.$%
\blacksquare $
\end{Beweis}

\subsection*{Continuous dual comodules}

\qquad In what follows we consider the \emph{dual comodules}{\normalsize \ }%
of modules of an $\alpha $-algebra over an arbitrary noetherian base ring.
These were considered in the case of base fields by several authors (e.g. 
\cite{Tak74}, \cite{Liu94}, \cite{GK97} and \cite{Lu98}) and in the case of
Dedekind domains by R. Larson \cite{Lar98}.

\begin{punto}
\label{0-funk}Let $R$ be noetherian, $A$ be an $R$-algebra, $\frak{F}\ $be a
filter consisting of $R$-cofinite $A$-ideals and consider $A$ with the right
induced linear topology $\frak{T}(\frak{F}).$

\begin{enumerate}
\item  If $\frak{F}\ $is an $\alpha $-filter basis, then by Proposition \ref
{A^0-F-co} (2) $A_{\frak{F}}^{\circ }$ is an $R$-coalgebra and $(A,A_{\frak{F%
}}^{\circ })\in \mathcal{P}_{m}^{\alpha }.$ By Theorem \ref{adj} we get
right-adjoint contravariant functors 
\begin{equation*}
\begin{tabular}{llllllll}
$(-)^{\ast }$ & $:$ & $\mathcal{M}^{A_{\frak{F}}^{\circ }}$ & $\rightarrow $
& $\mathcal{M}_{A},$ & $M$ & $\mapsto $ & $M^{\ast },$ \\ 
$(-)_{\frak{F}}^{0}$ & $:$ & $\mathcal{M}_{A}$ & $\rightarrow $ & $\mathcal{M%
}^{A_{\frak{F}}^{\circ }},$ & $M$ & $\mapsto $ & $M_{\frak{F}}^{0}:=\mathrm{%
Rat}^{A_{\frak{F}}^{\circ }}(_{A}M^{\ast }).$%
\end{tabular}
\end{equation*}
For every $M\in \mathcal{M}_{A}$ we call $M_{\frak{F}}^{0}$ the \emph{dual
comodule of }$M$\emph{\ w.r.t.}\textbf{\ }$\frak{F}$. If $A$ is a $\alpha $%
-algebra, then we call $M^{0}:=\mathrm{Rat}^{A^{\circ }}(_{A}M^{\ast })$ the 
\emph{dual comodule of }$M.$

\item  For every right $A$-module $M$ we call 
\begin{equation*}
M_{\frak{F}}^{\circ }\text{ }:=\{f\in M^{\ast }|\text{ }f(MI)=0\text{ for
some }I\in \frak{F}\}={\underrightarrow{lim}}_{\frak{F}}(M/MI)^{\ast }
\end{equation*}
the \emph{continuous dual module \ of }$M$\emph{\ w.r.t. }$\frak{F}\mathbf{.}
${\normalsize \ }If $A_{\frak{F}}^{\circ }$ is an $R$-coalgebra and $M_{%
\frak{F}}^{\circ }$ is a right $A_{\frak{F}}^{\circ }$-comodule, then we
call $M_{\frak{F}}^{\circ }$ the \emph{continuous dual comodule} of $M$
w.r.t. $\frak{F}.$
\end{enumerate}
\end{punto}

\begin{Notation}
Let $R$ be noetherian, $A$ be an ($\alpha $)-algebra and $M,N$ be right $A$%
-modules. For every $A$-linear mapping $\gamma :M\rightarrow N$ we denote
with $\gamma ^{\circ }:N^{\circ }\rightarrow M^{\circ }$ ($\gamma
^{0}:N^{0}\rightarrow M^{0}$) the restriction of $\gamma ^{\ast }$ on $%
N^{\circ }$ (on $N^{0}$).
\end{Notation}

\qquad The following result generalizes the corresponding one 
\cite[Corollary 2.2.16]{DNR} stated for the canonical pairing $(C^{\ast },C)$
over a base field to an arbitrary measuring $\alpha $-pairing $(A,C)$ over
an arbitrary noetherian ground ring:

\begin{proposition}
\label{N*^rat}Let $P=(A,C)\in \mathcal{P}_{m}^{\alpha },$ $N\in $ $^{C}%
\mathrm{Rat}(\mathcal{M}_{A})$ and consider for every $f\in N^{\ast }$ the $%
R $-linear mapping 
\begin{equation*}
\theta _{f}:N\rightarrow C,\text{ }\theta _{f}(n)=\sum n_{<-1>}f(n_{<0>}).
\end{equation*}
If $R$ is noetherian, then: 
\begin{equation*}
\begin{tabular}{lll}
$\mathrm{Rat}^{C}(_{A}N^{\ast })$ & $=$ & $\mathrm{Sp}(\sigma \lbrack
_{A}C], $ $_{A}N^{\ast }):=\sum \{\mathrm{Im}(g)|$ $g\in $ $\mathrm{Hom}%
_{A-}(U,N^{\ast }),$ $U\in \sigma \lbrack _{A}C]\}$ \\ 
& $=$ & $\{f\in N^{\ast }|$ $Af$ is finitely generated$\}$ \emph{(}$=$ $%
\mathrm{Loc}(_{A}N^{\ast })$\emph{)} \\ 
& $=$ & $\{f\in N^{\ast }|$ $\exists $ an $R$-cofinite \emph{(}right\emph{)}
ideal $I\subset A$ with $f(NI)=0\}$ \\ 
& $=$ & $\{f\in N^{\ast }|$ $\exists $ an $R$-cofinite $A$-submodule $%
L\subset N$ with $f(L)=0\}$ \\ 
& $=$ & $\{f\in N^{\ast }|$ $\exists $ an $R$-cofinite $C$-subcomodule $%
L\subset N$ with $f(L)=0\}$ \\ 
& $=$ & $\{f\in N^{\ast }|$ $\theta _{f}(N)\subset C$ is a finitely
generated $R$-submodule$\}.$%
\end{tabular}
\end{equation*}
\end{proposition}

\begin{Beweis}
The equality $\mathrm{Rat}^{C}(_{A}N^{\ast })=\mathrm{Sp}(\sigma \lbrack
_{A}C],$ $_{A}N^{\ast })$ follows from \ref{C-ad} and Theorem \ref{cor-dicht}%
. Obviously $\mathrm{Rat}^{C}(_{A}N^{\ast })\subseteq \mathrm{Loc}%
(_{A}N^{\ast }).$

By Theorem \ref{cor-dicht} and Lemma \ref{grun} $f\in \mathrm{Loc}%
(_{A}N^{\ast })\Leftrightarrow f(NI)=0$ for an $R$-cofinite (right) ideal $%
I\vartriangleleft A\Leftrightarrow f(L)=0$ for an $R$-cofinite right $A$%
-submodule $L\subset N$ $\Leftrightarrow $ $f(L)=0$ for an $R$-cofinite left 
$C$-subcomodule $L\subset N.$

Let $f\in N^{\ast }$ with $f(L)=0$ for an $R$-cofinite left $C$-subcomodule $%
L\subset N.$ Analogous to \ref{th-m} below, $\theta _{f}:N\rightarrow C$ is $%
C$-colinear. Notice that $\theta _{f}(L)=0$ and so there exists a $C$%
-colinear morphism $\overline{\theta _{f}}:N/L\rightarrow C,$ such that $%
\overline{\theta _{f}}\circ \pi _{L}=\theta _{f}.$ Consequently $\theta
_{f}(N)=\overline{\theta _{f}}(N/L)$ is finitely generated in $\mathcal{M}%
_{R}.$

To every $f\in N^{\ast }$ there corresponds the left $C$-coideal $\theta
_{f}(N)\subset C.$ If $\theta _{f}(N)$ is finitely generated in $\mathcal{M}%
_{R},$ then $(\theta _{f}(N))^{\ast }$ is a right $C$-comodule by \ref{fp}
and we have for every $n\in N:$%
\begin{equation*}
\varepsilon (\theta _{f}(n))=\varepsilon (\sum f(n_{<0>})n_{<-1>})=f(\sum
\varepsilon (n_{<-1>})n_{<0>})=f(n),
\end{equation*}
\newline
i.e. $f\in (\theta _{f}(N))^{\ast }\subset \mathrm{Rat}^{C}(_{A}N^{\ast
}).\blacksquare $
\end{Beweis}

\qquad As a special case of Proposition \ref{N*^rat} we get

\begin{corollary}
\label{C*^rat}Let $R$ be noetherian. For every locally projective $R$%
-coalgebra $C$ we have 
\begin{equation*}
\begin{tabular}{lll}
$\mathrm{Rat}^{C}(_{C^{\ast }}C^{\ast })$ & $=$ & $\mathrm{Sp}(\sigma
\lbrack _{C^{\ast }}C],$ $_{C^{\ast }}C^{\ast }):=\sum \{\mathrm{Im}(g)|$ $%
g\in $ $\mathrm{Hom}_{C^{\ast }-}(U,C^{\ast }),$ $U\in \sigma \lbrack
_{C^{\ast }}C]\}$ \\ 
& $=$ & $\{f\in C^{\ast }|$ $C^{\ast }\star f$ is finitely generated in $%
\mathcal{M}_{R}\}$ \\ 
& $=$ & $\{f\in C^{\ast }|$ $\exists $ an $R$-cofinite \emph{(}right\emph{)}
ideal $I\vartriangleleft C^{\ast }$ with $f(CI)=0\}$ \\ 
& $=$ & $\{f\in C^{\ast }|$ $\exists $ an $R$-cofinite right $C^{\ast }$%
-submodule $K\subset C$ with $f(K)=0\}$ \\ 
& $=$ & $\{f\in C^{\ast }|$ $\exists $ an $R$-cofinite left $C$-coideal $%
K\subset C$ with $f(K)=0\}$ \\ 
& $=$ & $\{f\in C^{\ast }|$ $f\rightharpoonup C\subset C$ is a finitely
generated $R$-submodule$\}.$%
\end{tabular}
\end{equation*}
\end{corollary}

\begin{punto}
\textbf{Cofree comodules.}\label{ot C-adj} A right $C$-comodule $(M,\varrho
_{M})$ is called \emph{cofree}, if there exists an $R$-module $K,$ such that 
$(M,\varrho _{M})\simeq (K\otimes _{R}C,id_{K}\otimes \Delta _{C})$ as right 
$C$-comodules. Note that if $K=R^{(\Lambda )},$ a free $R$-module, then $%
M\simeq R^{(\Lambda )}\otimes _{R}C\simeq C^{(\Lambda )}$ as right $C$%
-comodules (in fact, this is the reason for the terminology $\emph{cofree}$).
\end{punto}

\begin{lemma}
\label{co-free}Let $R$ be noetherian and $A$ be a cofinitary\emph{\ }$R$%
-algebra. Let $M$ be an $R$-module and consider the right $A$-module $%
N:=M\otimes _{R}A.$ Then $N^{\circ }\simeq M^{\ast }\otimes _{R}A^{\circ }$
as $A^{\circ }$-comodules \emph{(}i.e. $N^{\circ }$ is a \emph{cofree}%
{\normalsize \ }right $A^{\circ }$-comodule\emph{)}.
\end{lemma}

\begin{Beweis}
If $N\simeq M\otimes _{R}A$ as right $A$-modules, then there are
isomorphisms in $_{A}\mathcal{M}:$%
\begin{equation*}
\begin{tabular}{llll}
$N^{\circ }$ & $:=$ & $\underrightarrow{lim}\{((M\otimes _{R}A)/(M\otimes
_{R}A)I)^{\ast }|$ $I\in \mathcal{K}_{A}\}$ &  \\ 
& $=$ & $\underrightarrow{lim}\{((M\otimes _{R}A)/(M\otimes _{R}A)\widetilde{%
I})^{\ast }|$ $\widetilde{I}\in \mathcal{E}_{A}\}$ & ($A$ is cofinitary) \\ 
& $=$ & $\underrightarrow{lim}\{((M\otimes _{R}A)/(M\otimes _{R}\widetilde{I}%
))^{\ast }|$ $\widetilde{I}\in \mathcal{E}_{A}\}$ &  \\ 
& $\simeq $ & $\underrightarrow{lim}\{(M\otimes _{R}A/\widetilde{I})^{\ast
}| $ $\widetilde{I}\in \mathcal{E}_{A}\}$ &  \\ 
& $\simeq $ & $\underrightarrow{lim}\{M^{\ast }\otimes _{R}(A/\widetilde{I}%
)^{\ast }|$ $\widetilde{I}\in \mathcal{E}_{A}\}$ & ($A/\widetilde{I}$ is
f.g. projective in $\mathcal{M}_{R}$); \\ 
& $\simeq $ & $M^{\ast }\otimes _{R}\underrightarrow{lim}\{(A/\widetilde{I}%
)^{\ast }|$ $\widetilde{I}\in \mathcal{E}_{A}\}$ &  \\ 
& $\simeq $ & $M^{\ast }\otimes _{R}A^{\circ }$ & ($A$ is cofinitary).$%
\blacksquare $%
\end{tabular}
\end{equation*}
\end{Beweis}

\qquad In contradiction with \cite[Corollary 2]{Wit79} the following example
shows that for an arbitrary $R$-algebra $A$ the preradical $\mathrm{Loc}(-):$
$_{A}\mathcal{M}\rightarrow $ $\mathrm{Loc}(_{A}\mathcal{M})$ is in general
not a torsion radical:

\begin{c-example}
\label{Rxft}(Compare \cite[Seite 155]{Mon93}) Let $R$ be a field and
consider the Hopf $R$-algebra $H:=R[x_{1},x_{2},...,x_{n},...],$ with the
usual multiplication, the usual unity and the comultiplication, counity and
antipode defined on the generators through 
\begin{equation*}
\Delta (x_{i})=1\otimes x_{i}+x_{i}\otimes 1,\text{ }\varepsilon (x_{i})=0,%
\text{ }S(x^{i}):=(-1)^{i}x^{i}.
\end{equation*}
If we consider $H$ with the left cofinite topology, then $(H,\mathrm{Cf}(H))$
is a left linear topological $R$-algebra with preradical $\mathrm{Loc}(-):$ $%
_{H}\mathcal{M}\rightarrow $ $_{H}\mathcal{M}$ and pretorsion class $\mathrm{%
Loc}(_{H}\mathcal{M})$ (see \ref{ko-top} and Proposition \ref{A^0-cor}). If
we consider the $H$-ideal $\omega :=Ke(\varepsilon _{H}),$ then $H/\omega
\simeq R$ while $\mathrm{dim}(H/\omega ^{2})=\infty ,$ i.e. $\omega
^{2}\notin \mathcal{K}_{H}.$ So $\mathrm{Cf}(H)$ is not a Gabriel-topology
and consequently $\mathrm{Loc}(_{H}\mathcal{M})$ is not closed against
extensions (see \cite[Chapter VI, Theorem 5.1, Lemma 5.3]{Ste75}).
\end{c-example}

\section{Coreflexive comodules}

\qquad In \cite{Taf72}, \cite{Taf77} E. Taft developed an algebraic aspect
to the study of \emph{coreflexive coalgebras} over base field (i.e.
coalgebras $C$ with $C\simeq C^{\ast \circ }$). Independently, R. Heyneman
and D. Radford \cite{Rad73}, \cite{HR74} studied the coreflexive coalgebras
with the help of the \emph{finite topology}\ on $C^{\ast }.$ In this section
we present and study for every $(A,C)\in \mathcal{P}_{m}^{\alpha }$ over an
arbitrary noetherian ring the notions of \emph{reflexive }$A$\emph{-modules}
and \emph{coreflexive\ }$C$\emph{-comodules. }We get algebraic as well as
topological characterizations for the (co)reflexive (co)modules. Our results
will be applied then to the study of (co)reflexive (co)algebras, where we
generalize also results from the papers mentioned above and from \cite{Wit79}%
.

\subsection*{$(A,C)$-pairings}

In the case of base fields, D. Radford \cite{Rad73} presented for every
measuring $R$-pairing $P=(A,C)\;$the so called right (left) $P$-pairings. In
what follows we consider \emph{duality relations} for such pairings.

\begin{punto}
Let $P=(A,C)\in \mathcal{P}_{m}.$ A pairing of $R$-modules $Q=(M,N)$ is
called a \emph{right }(\emph{a left})\emph{\ }$P$\emph{-pairing}, if $M$ is
a right (a left) $A$-module, $N$ is a right (a left) $C$-comodule and the
induced mapping $\kappa _{Q}:M\rightarrow N^{\ast }$ is $A$-linear. By $%
\mathcal{Q}_{P}^{r}\subset \mathcal{P}$ ($\mathcal{Q}_{P}^{l}\subset 
\mathcal{P}$) we denote the subcategory of right (left) $P$-pairings with
morphisms described as follows: if $(M,N),$ $(M^{\prime },N^{\prime })$ are
right (left) $P$-pairings, then a morphism of $R$-pairings 
\begin{equation*}
(\xi ,\theta ):(M^{\prime },N^{\prime })\rightarrow (M,N)
\end{equation*}
is a morphism in $\mathcal{Q}_{P}^{r}$ (resp. in $\mathcal{Q}_{P}^{l}$), if $%
\xi :M\rightarrow M^{\prime }$ is $A$-Linear and $\theta :N^{\prime
}\rightarrow N$ is $C$-colinear.

A $P$\emph{-bi-pairing} is an $R$-pairing $(M,N),$ where $M$ is an $A$%
-bimodule, $N$ is a $C$-bicomodule and $\kappa _{Q}:M\rightarrow N^{\ast }$
is $A$-bilinear. With $\mathcal{Q}_{P}$ we denote the category of $P$%
-bi-pairings with morphisms described as follows: if $(M,N),$ $(M^{\prime
},N^{\prime })$ are $P$-bi-pairings, then a morphism of $R$-pairings 
\begin{equation*}
(\xi ,\theta ):(M^{\prime },N^{\prime })\rightarrow (M,N)
\end{equation*}
is a morphism in $\mathcal{Q}_{P},$ if $\xi :M\rightarrow M^{\prime }$ is $A$%
-bilinear and $\theta :N^{\prime }\rightarrow N$ is $C$-bicolinear. In
particular every measuring $R$-pairing $P$ is itself a $P$-bi-pairing.
\end{punto}

\begin{punto}
\label{th-m}Let $P=(A,C)\in \mathcal{P}_{m},$ $Q=(M,N)\in \mathcal{Q}%
_{P}^{r} $ and define for every $m\in M:$ 
\begin{equation*}
\begin{tabular}{llllllll}
$\xi _{m}:$ & $A$ & $\rightarrow $ & $M,$ & $a$ & $\mapsto $ & $ma$ & $\text{
for all }a\in A,$ \\ 
$\theta _{m}:$ & $N$ & $\rightarrow $ & $C,$ & $n$ & $\mapsto $ & $\sum
<m,n_{<0>}>n_{<1>}$ & $\text{ for all }n\in N.$%
\end{tabular}
\end{equation*}
Then we have for all $a\in A$ and $n\in N:$%
\begin{equation*}
<\xi _{m}(a),n>=<ma,n>=<m,an>=\sum <m,n_{<0>}><a,n_{<1>}>=<a,\theta _{m}(n)>.
\end{equation*}
Obviously $\xi _{m}:A\rightarrow M$ is $A$-linear. Moreover it follows for
all $n\in N$ and $a\in A$ that 
\begin{equation*}
\begin{tabular}{lll}
$\alpha _{N}^{P}(\sum \theta _{m}(n)_{1}\otimes \theta _{m}(n)_{2})(a)$ & $=$
& $\sum \theta _{m}(n)_{1}<a,\theta _{m}(n)_{2}>$ \\ 
& $=$ & $a\rightharpoonup \theta _{m}(n)$ \\ 
& $=$ & $\sum <m,n_{<0>}>(a\rightharpoonup n_{<1>})$ \\ 
& $=$ & $\sum <m,n_{<0>}>n_{<1>1}<a,n_{<1>2}>$ \\ 
& $=$ & $\sum <m,n_{<0><0>}>n_{<0><1>}<a,n_{<1>}>$ \\ 
& $=$ & $\alpha _{N}^{P}(\sum \theta _{m}(n_{<0>})\otimes n_{<1>})(a).$%
\end{tabular}
\end{equation*}
If $\alpha _{N}^{P}:N\otimes _{R}C\rightarrow \mathrm{Hom}_{R}(A,N)$ is
injective, then 
\begin{equation*}
\sum \theta _{m}(n)_{1}\otimes \theta _{m}(n)_{2}=\sum \theta
_{m}(n_{<0>})\otimes n_{<1>}\text{ for every }n\in N,
\end{equation*}
i.e. $\theta _{m}:N\rightarrow C$ is $C$-colinear and 
\begin{equation*}
(\xi _{m},\theta _{m}):(M,N)\rightarrow (A,C)
\end{equation*}
is a morphism in $\mathcal{Q}_{P}^{r}.$
\end{punto}

\begin{Notation}
Let $P=(A,C)\in \mathcal{P}_{m}$ and $Q=(M,N)\in \mathcal{Q}_{P}^{r}.$ For $%
R $-submodules $L\subset M,$ $K\subset N$ we set 
\begin{equation*}
\begin{tabular}{llll}
$K^{\bot }:=$ & $\{m\in M|$ $<m,K>=0\},$ & $\mathrm{An}_{M}(K):=$ & $\{m\in
M|\text{ }\theta _{m}(k)=0$ $\forall $ $k\in K\},$ \\ 
$L^{\bot }:=$ & $\{n\in N|$ $<L,n>=0\},$ & $\mathrm{An}_{N}(L):=$ & $\{n\in
N|$ $\theta _{m}(n)=0$ $\forall $ $m\in L\}.$%
\end{tabular}
\end{equation*}
\end{Notation}

\qquad As a consequence of Theorem \ref{cor-dicht} one can easily derive the
following result:

\begin{lemma}
\label{mod-per}Let $P=(A,C)\in \mathcal{P}_{m}$ and $Q=(M,N)\in \mathcal{Q}%
_{P}^{r}$ \emph{(}resp. $Q\in \mathcal{Q}_{P}^{l},$ $Q\in \mathcal{Q}_{P}$%
\emph{)}.

\begin{enumerate}
\item  Every right $C$-subcomodule \emph{(}resp. left $C$-subcomodule, $C$%
-subbicomodule\emph{)} $K\subset N$ is a left $A$-submodule \emph{(}resp. a
right $A$-submodule, an $A$-subbimodule\emph{)}\ and ${K}^{\bot }\subset M$
is a right $A$-submodule \emph{(}resp. a left $A$-submodule, an $A$%
-subbimodule\emph{)}.

\item  Let $(A,C)\in \mathcal{P}_{m}^{\alpha }.$ If $L\subset M$ a right $A$%
-submodule \emph{(}resp. a left $A$-submodule, an $A$-subbimodule\emph{)},
then ${L}^{\bot }\subset N$ is a right $C$-subcomodule \emph{(}resp. a left $%
C$-subcomodule, a $C$-subbicomodule\emph{)}.
\end{enumerate}
\end{lemma}

\subsection*{{The topology $\frak{T}_{N}(M)$}}

\qquad

Let $P=(A,C)$ be a measuring $R$-pairing and consider $A$ as a right linear
topological $R$-algebra with the right $C$-adic topology $\mathcal{T}%
_{-C}(A).$ For every $Q=(M,N)\in \mathcal{Q}_{P}^{r}$ we present on $M$ a
topology $\frak{T}_{N}(M),$ such that $(M,\frak{T}_{N}(M))$ is a linear
topological $(A,\mathcal{T}_{-C}(A))$-module.

\begin{punto}
Let $P=(A,C)\in \mathcal{P}_{m},$ $Q=(M,N)\in \mathcal{Q}_{P}^{r}$ and
consider $C$ with the canonical right $A$-module structure and $A$ as a
right linear topological $R$-algebra with the right $C$-adic topology $%
\mathcal{T}_{-C}(A)$ (compare \ref{C-ad}). If $K\subset N$ is an $R$%
-submodule and $m\in \mathrm{An}_{M}(K),$ then we have for arbitrary $n\in K$
and $a\in A:$%
\begin{equation*}
\begin{tabular}{lllll}
$\theta _{ma}(n)$ & $:=$ & $\sum <ma,n_{<0>}>n_{<1>}$ &  &  \\ 
& $=$ & $\sum <m,an_{<0>}>n_{<1>}$ &  &  \\ 
& $=$ & $\sum <m,n_{<0><0>}><a,n_{<0><1>}>n_{<1>}$ &  &  \\ 
& $=$ & $\sum <m,n_{<0>}><a,n_{<1>1}>n_{<1>2}$ &  &  \\ 
& $=$ & $[(\alpha _{C}^{P}\circ \Delta _{C}^{cop})(\sum
<m,n_{<0>}>n_{<1>})](a)$ &  &  \\ 
& $=$ & $[(\alpha _{C}^{P}\circ \Delta _{C}^{cop})(\theta _{m}(n))](a)=0,$ & 
& 
\end{tabular}
\end{equation*}
i.e. $\mathrm{An}_{M}(K)\subset M$ is an $A$-submodule. Let $%
K=\sum\limits_{i=1}^{l}Rn_{i}\subset N$ be an arbitrary finitely generated $%
R $-submodule with $\varrho
_{N}(n_{i})=\sum\limits_{j=1}^{l_{i}}n_{ij}\otimes c_{ij}$ for $i=1,...,l$
and set $W:=\sum\limits_{i=1}^{l}\sum\limits_{j=1}^{l_{i}}Rc_{ij}.$ Let $%
m\in M$ be arbitrary. If $a\in \mathrm{An}_{A}^{r}(W),$ then for $i=1,...,l:$%
\begin{equation*}
\begin{tabular}{lllll}
$\theta _{ma}(n_{i})$ & $=$ & $\sum\limits_{i=1}^{l}\sum%
\limits_{j=1}^{l_{i}}<ma,n_{ij}>c_{ij}$ &  &  \\ 
& $=$ & $\sum\limits_{i=1}^{l}\sum\limits_{j=1}^{l_{i}}<m,an_{ij}>c_{ij}$ & 
&  \\ 
& $=$ & $\sum\limits_{i=1}^{l}\sum\limits_{j=1}^{l_{i}}\sum%
\limits_{n_{ij}}<m,n_{ij<0>}><a,n_{ij<1>}>c_{ij}$ &  &  \\ 
& $=$ & $\sum\limits_{i=1}^{l}\sum\limits_{j=1}^{l_{i}}\sum%
\limits_{c_{ij}}<m,n_{ij}><a,c_{ij1}>c_{ij2}$ &  &  \\ 
& $=$ & $\sum\limits_{i=1}^{l}\sum\limits_{j=1}^{l_{i}}<m,n_{ij}>(c_{ij}%
\leftharpoonup a)=0,$ &  & 
\end{tabular}
\end{equation*}
i.e. $(\mathrm{An}_{M}(K):m)\supset \mathrm{An}_{A}^{r}(W)$ and so it's open
w.r.t. the $C$-adic topology $\mathcal{T}_{-C}(A).$ So 
\begin{equation*}
\mathcal{B}(0_{M}){\normalsize :=\{\mathrm{An}_{M}(K)|}\text{ }K\subset N%
\text{ is a finitely generated }R\text{-submodule}\}
\end{equation*}
is neighbourhood basis of $0_{M}$ consisting of $A$-submodules of $M$ and $M$
becomes a topology $\frak{T}_{N}(M),$ such that $(M,\frak{T}_{N}(M))$ is a
linear topological right $(A,\mathcal{T}_{-C}(A))$-module.
\end{punto}

\begin{punto}
Let $P=(A,C)\in \mathcal{P}_{m}$ and $Q=(M,N)\in \mathcal{Q}_{P}^{r}.$ Then 
\begin{equation*}
\mathcal{F}{\normalsize (0_{M}):=\{K^{\bot }|}\text{ }K\subset N\text{ is a
finitely generated }R\text{-submodule}\}
\end{equation*}
is a filter basis consisting of $R$-submodules of $M$ and induces on $M$ the 
\emph{linear weak topology} $M[\frak{T}_{ls}(N)],$ such that $(M,M[\frak{T}%
_{ls}(N)])$ is a linear topological $R$-module and $\mathcal{F}(0_{M})$ is a
neighbourhood basis of $0_{M}.$
\end{punto}

\begin{punto}
Let $R$ be noetherian, $P=(A,C)\in \mathcal{P}_{m}$ and consider $C^{\ast }$
with the \emph{right cofinite topology} $\mathrm{Cf}(C^{\ast })$ (see \ref
{ko-top}). The $R$-algebra morphism $\kappa _{P}:A\rightarrow C^{\ast }$
induces on $A$ a linear topology $\kappa _{P}$-$\mathrm{Cf}(A)$ with
neighbourhood basis of $0_{A}:$%
\begin{equation*}
\mathcal{B}_{\kappa _{P}}(0_{A}):=\{\kappa _{P}^{-1}(J)|\text{ }%
J\vartriangleleft C^{\ast }\text{ is an $R$-cofinite two-sided ideal}\}.
\end{equation*}
By definition $\kappa _{P}$-$\mathrm{Cf}(A)$ the \emph{finest} linear
topology $\frak{T}$ on $A,$ such that $(A,\frak{T})$ is a linear topological 
$R$-algebra and $\kappa _{P}:A\rightarrow C^{\ast }$ is continuous.

Let $Q=(M,N)\in \mathcal{Q}_{P}^{r}$ and consider $N_{A}^{\ast }$ with the
cofinite topology $\mathrm{Cf}(N^{\ast }).$ The $A$-linear mapping $\kappa
_{Q}:M\rightarrow N^{\ast }$ induces on $M$ a topology $\kappa _{Q}$%
{\normalsize -}$\mathrm{Cf}(M)$ with neighbourhood basis of $0_{M}$%
\begin{equation*}
\mathcal{B}_{\kappa _{Q}}(0_{M}):=\{\kappa _{Q}^{-1}(L)|\text{ }L\subset
N^{\ast }\text{ is an $R$-cofinite }A\text{-submodule}\}.
\end{equation*}
Clearly $\kappa _{Q}$-$\mathrm{Cf}(M)$ is a linear topological right $%
\mathrm{Cf}(A)$-module and is the \emph{finest} topology $\frak{T}$ on $M,$
such that $(M,\frak{T})$ is a linear topological right $(A,\mathrm{Cf}(A))$%
-module and $\kappa _{Q}:M\rightarrow N^{\ast }$ is continuous.
\end{punto}

\begin{lemma}
\label{leq}Let $P=(A,C)\in \mathcal{P}_{m}$ and $Q=(M,N)\in \mathcal{Q}%
_{P}^{r}.$

\begin{enumerate}
\item  The linear weak topology $M[\frak{T}_{ls}(N)]$ and the topology $%
\frak{T}_{N}(M)$ coincide. So $M,$ considered with the linear weak topology,
is a linear topological right $(A,\mathcal{T}_{-C}(A))$-module.

\item  If $R$ is noetherian and $P$ satisfies the $\alpha $-condition, then 
\begin{equation}
M[\frak{T}_{ls}(N)]\leq \kappa _{Q}\text{-}\mathrm{Cf}(M)\leq \mathrm{Cf}(M).
\label{top-rel}
\end{equation}
\end{enumerate}
\end{lemma}

\begin{Beweis}
\begin{enumerate}
\item  Let $U\subset M$ be a neighbourhood of $0_{M}$ w.r.t. $M[\frak{T}%
_{ls}(N)].$ Then there exists a finitely generated $R$-submodule $K\subset
N, $ such that $K^{\bot }\subseteq U.$ If $m\in \mathrm{An}_{M}(K),$ then we
have for arbitrary $n\in K:$%
\begin{equation*}
<m,n>=<m,\sum n_{<0>}\varepsilon (n_{<1>})>=\varepsilon (\sum
<m,n_{<0>}>n_{<1>})=\varepsilon (\theta _{m}(n))=0.
\end{equation*}
So $\mathrm{An}_{M}(K)\subseteq K^{\bot }\subseteq U,$ i.e. $U$ is a
neighbourhood of $0_{M}$ w.r.t. $\frak{T}_{N}(M).$

On the other hand, let $U\subset M$ be a neighbourhood of $0_{M}$ w.r.t. $%
\frak{T}_{N}(M).$ Then there exists a finitely generated $R$-submodule $%
K=\sum\limits_{i=1}^{l}Rn_{i}\subset N,$ such that $\mathrm{An}%
_{M}(K)\subseteq U.$ Assume now that $\varrho
_{N}(n_{i})=\sum\limits_{j=1}^{l_{i}}n_{ij}\otimes c_{ij}$ and set $%
W:=\sum\limits_{i=1}^{l}\sum\limits_{j=1}^{l_{i}}Rn_{ij}.$ Then $W^{\bot
}\subseteq \mathrm{An}_{M}(K)\subseteq U,$ i.e. $U$ is a neighbourhood of $%
0_{M}$ w.r.t. $M[\frak{T}_{ls}(N)].$ Consequently $M[\frak{T}_{ls}(N)]=\frak{%
T}_{N}(M).$

\item  Let $R$ be noetherian and $P\in \mathcal{P}_{m}^{\alpha }.$ Let $%
U\subset M$ be a neighbourhood of $0_{M}$ w.r.t. $M[\frak{T}_{ls}(N)],$ i.e.
there exists a finitely generated $R$-submodule $K\subset N,$ such that $%
K^{\bot }\subseteq U.$ By assumption $P\in \mathcal{P}_{m}^{\alpha }$ and so
there exists by the \emph{Finiteness Theorem} \ref{es} a left $A$-submodule $%
\widetilde{K}\subset N,$ such that $K\subseteq \widetilde{K}$ and $%
\widetilde{K}_{R}$ is finitely generated. Moreover $N^{\ast }/\mathrm{An}(%
\widetilde{K})\hookrightarrow \widetilde{K}^{\ast },$ i.e. $\mathrm{An}(%
\widetilde{K})\subset N^{\ast }$ is an $R$-cofinite right $A$-submodule and $%
\kappa _{Q}^{-1}(\mathrm{An}(\widetilde{K})):=\widetilde{K}^{\bot }\subseteq
K^{\bot }\subset U,$ i.e. $U$ is a neighbourhood of $0_{M}$ w.r.t. $\kappa
_{Q}$-$\mathrm{Cf}(M).$

Let $U\subset M$ be a neighbourhood of $0_{M}$ w.r.t. $\kappa _{Q}$-$\mathrm{%
Cf}(M),$ i.e. there exists an $R$-cofinite $A$-submodule $L\subset N^{\ast
}, $ such that $\kappa _{Q}^{-1}(L)\subseteq U.$ Then $M/\kappa
_{Q}^{-1}(L)\hookrightarrow N^{\ast }/L,$ and so $\kappa _{Q}^{-1}(L)\subset
M$ is an $R$-cofinite $A$-submodule. Consequently $U$ is a neighbourhood of $%
0_{M}$ w.r.t. $\mathrm{Cf}(M).\blacksquare $
\end{enumerate}
\end{Beweis}

\begin{definition}
Let $P=(A,C)\in \mathcal{P}_{m}$ and $Q=(M,N)\in \mathcal{Q}_{P}^{r}.$

\begin{enumerate}
\item  If $P\in \mathcal{P}_{m}^{\alpha },$ then we call $Q$ \emph{weakly
coreflexive}, if $N=M^{r}.$

\item  If $R$ is noetherian, then we call $Q$ \emph{coreflexive}, if $M[%
\frak{T}_{ls}(N)]=\mathrm{Cf}(M).$

\item  We call $Q$ \emph{proper} (resp. \emph{weakly reflexive, reflexive}),
if $\kappa _{Q}:M\rightarrow N^{\ast }$ is injective (resp. surjective,
bijective).
\end{enumerate}
\end{definition}

\begin{definition}
\begin{enumerate}
\item  Let $C$ be an $R$-coalgebra and $N$ be a right $C$-comodule.

\begin{enumerate}
\item  If $_{R}C$ is locally projective, then we call $N$ \emph{weakly
coreflexive,} if $N=N^{\ast \Box }.$

\item  If $R$ is noetherian, then we call $N$ \emph{coreflexive}, if $%
N^{\ast }[\frak{T}_{ls}(N)]=\mathrm{Cf}(N^{\ast }).$
\end{enumerate}

\item  Let $R$ be noetherian and $A$ be an $R$-algebra. We call a right $A$%
-module $M$ \emph{proper} (resp. \emph{weakly reflexive, reflexive}), if the
canonical $A$-linear mapping $\lambda _{M}:M\rightarrow M^{\circ \ast }$ is
injective (resp. surjective, bijective).
\end{enumerate}
\end{definition}

\begin{remarks}
\label{(A,C)-sr}

\begin{enumerate}
\item  Consider the ground ring $R$ as a trivial $R$-bialgebra. Then $%
R^{\ast }\simeq R,$ $\mathcal{M}_{R}\simeq \mathcal{M}^{R}$ and for every $R$%
-(co-)module $N$ we have $N^{\ast \ast }=\mathrm{Rat}^{R}(N^{\ast \ast })=%
\mathrm{Loc}(_{R^{\ast }}N^{\ast \ast }).$ So $N$ is (co)reflexive, iff $N$
is reflexive in the usual sense, i.e. if the canonical $R$-linear mapping $%
\Phi _{N}:N\rightarrow N^{\ast \ast }$ is bijective.

\item  For every $P=(A,C)\in \mathcal{P}_{m}^{\alpha }$ we have $C=A^{r}$
(by Corollary \ref{C=A^r} (1)) and so $P\in \mathcal{Q}_{P}^{r}$ is weakly
coreflexive.
\end{enumerate}
\end{remarks}

\begin{proposition}
\label{A-prop}Let $R$ be noetherian and $A$ be an $R$-algebra.

\begin{enumerate}
\item  If $A$ is proper \emph{(}i.e. the canonical mapping $\lambda
_{A}:A\rightarrow A^{\circ \ast }$ is injective\emph{)}, then $\mathrm{Cf}%
(A) $ is Hausdorff.

\item  Let $A$ be cofinitely $R$-cogenerated. Then $A$ is proper, iff $%
\mathrm{Cf}(A)$ is Hausdorff.

\item  If $R$ is a QF ring, then 
\begin{equation*}
A\text{ is proper }\Leftrightarrow \text{ }\mathrm{Cf}(A)\text{ is Hausdorff}%
\Leftrightarrow \text{ }A^{\circ }\subset A^{\ast }\text{ is dense.}
\end{equation*}
\end{enumerate}
\end{proposition}

\begin{Beweis}
\begin{enumerate}
\item  Obviously $\overline{0_{A}}:=\bigcap\limits_{\mathcal{K}_{A}}I\subset 
\mathrm{Ke}(\lambda _{A})$ and the result follows.

\item  Assume $\mathrm{Cf}(A)$ to be Hausdorff. If $A$ is not proper, then
there exists some $0\neq \widetilde{a}\in A$, such that $f(\widetilde{a})=0$
for every $f\in A^{\circ }.$ If $I\vartriangleleft A$ is an arbitrary $R$%
-cofinite $A$-ideal, then $\widetilde{a}\in \mathrm{KeAn}(I)=I$ (compare 
\cite[28.1]{Wis88}) and so $\bigcap\limits_{\mathcal{K}_{A}}I\neq 0$ (\emph{%
contradiction}).

\item  By \cite[Theorem 1.8 (1)]{Abu} we have 
\begin{equation*}
\overline{A^{\circ }}=\mathrm{AnKe}(A^{\circ })=\mathrm{An}{\normalsize (}%
\mathrm{Ke}(\sum_{I\in \mathcal{K}_{A}}\mathrm{An}(I))=\mathrm{An}%
(\bigcap_{I\in \mathcal{K}_{A}}\mathrm{KeAn}(I))=\mathrm{An}(\bigcap_{I\in 
\mathcal{K}_{A}}I).
\end{equation*}
So $A^{\circ }\subset A^{\ast }$ is dense, iff $\bigcap\limits_{\mathcal{K}%
_{A}}I=0.\blacksquare $
\end{enumerate}
\end{Beweis}

\begin{lemma}
\label{krull}\emph{(}\textbf{Krull's Theorem}\emph{)} Let $A$ be a \emph{%
commutative noetherian }ring. For every \emph{finitely generated }$A$-module 
$M$ and every $A$-ideal $I\vartriangleleft A$ we have 
\begin{equation*}
\bigcap_{k=0}^{\infty }MI^{k+1}=\{m\in M|\text{ }\exists \text{ }b\in I,%
\text{ such that }m(1_{A}-b)=0\}.
\end{equation*}
\end{lemma}

\qquad The following result was obtained in \cite[6.1.3]{Swe69} for \emph{%
commutative affine algebras} over base fields:

\begin{lemma}
\label{class-prop}Let $R$ be a QF ring and $A$ be a \emph{commutative
noetherian} $R$-algebra. If every maximal $A$-ideal is $R$-cofinite, then $%
A^{\circ }\subset A^{\ast }$ is dense.
\end{lemma}

\begin{Beweis}
Let $0\neq a\in A$ be arbitrary and consider the $A$-ideal $J:=(0:a).$ Let $%
\frak{m}\vartriangleleft A$ be a maximal $A$-ideal, such that $J\subset 
\frak{m}.$ Sine $A$ is noetherian, $\frak{m}_{A}$ is finitely generated. If $%
a\in \bigcap\limits_{k=0}^{\infty }\frak{m}^{k+1},$ then there exists by
Krull's Theorem some $b\in \frak{m},$ such that $a(1_{A}-b)=0$ and so $%
1_{A}\in \frak{m}$ (\emph{contradiction}). So there exists $k\geq 0,$ such
that $a\notin \frak{m}^{k+1}.$ By assumption\ $\frak{m}\subset A$ is $R$%
-cofinite and it follows then from Lemma \ref{jm} that $\frak{m}%
^{k+1}\subseteq A$ is $R$-cofinite, i.e. $a\notin \bigcap_{\mathcal{K}%
_{A}}I. $ Since $0\neq a\in A$ is arbitrary by our choice, it follows that $%
\bigcap_{\mathcal{K}(A)}I=0,$ i.e. $A$ is proper and consequently $A^{\circ
}\subset A^{\ast }$ is dense by Proposition \ref{A-prop}.$\blacksquare $
\end{Beweis}

Analog to the proof of Proposition \ref{A-prop} we get

\begin{proposition}
\label{ls-prop}Let $R$ be noetherian, $A$ be an $R$-algebra and $M$ be a
right $A$-module.

\begin{enumerate}
\item  If $M$ is proper, then $\mathrm{Cf}(M)$ is Hausdorff.

\item  Let $M$ be cofinitely $R$-cogenerated. Then $\mathrm{Cf}(M)$ is
Hausdorff, iff $M$ is proper.

\item  If $R$ is a QF ring, then 
\begin{equation*}
M\text{ is proper }\Leftrightarrow \text{ }\mathrm{Cf}(M)\text{ is Hausdorff}%
\Leftrightarrow \text{ }M^{\circ }\subset M^{\ast }\text{ is dense.}
\end{equation*}
\end{enumerate}
\end{proposition}

\begin{theorem}
\label{Q-koref}Let $R$ be noetherian, $P=(A,C)\in \mathcal{P}_{m}^{\alpha }$
and $Q=(M,N)\in \mathcal{Q}_{P}^{r}.$

\begin{enumerate}
\item  If $Q$ is coreflexive, then $M^{r}=M^{\circ }.$

\item  Let $M$ be cofinitely $R$-cogenerated.

\begin{enumerate}
\item  If $N\overset{\chi _{Q}}{\simeq }M^{\circ }$, then $Q$ is coreflexive.

\item  Let $Q$ be weakly coreflexive. Then $Q$ is coreflexive, iff $N%
\overset{\chi _{Q}}{\simeq }M^{\circ }.$
\end{enumerate}
\end{enumerate}
\end{theorem}

\begin{Beweis}
\begin{enumerate}
\item  Assume $Q$ to be coreflexive and consider $A$ and $M$ with the linear
weak topology $A[\frak{T}_{ls}(C)],$ $M[\frak{T}_{ls}[N])$ respectively. Let 
$f\in M^{\ast }$ with $f(L)=0$ for an $R$-cofinite $A$-submodule $L\subset
M, $ say $M/L=\sum\limits_{i=1}^{k}R(m_{i}+L).$ By assumption $M[\frak{T}%
_{ls}(N)]=\mathrm{Cf}(M)$ and so $L$ is open w.r.t. $M[\frak{T}_{ls}(N)].$
By \cite[Corollary 1.9]{Abu} $\xi _{m_{i}}:A\rightarrow M$ is continuous and
so there exist finitely generated $R$-submodules $Z_{1},...,Z_{k}\subseteq
C, $ such that $Z_{i}^{\bot }\subseteq \xi _{m_{i}}^{-1}(L).$ Consequently $%
(\sum\limits_{i=1}^{k}Z_{i})^{\bot }=\bigcap\limits_{i=1}^{k}Z_{i}^{\bot
}\subseteq (0_{M^{\ast }}:f),$ i.e. $f\in M^{r}$ (by Proposition \ref{n_rat}%
). Obviously $M^{r}\subset M^{\circ }$ and the result follows.

\item  Let $M\ $be cofinitely $R$-cogenerated.

\begin{enumerate}
\item  Assume that $N\overset{\chi _{Q}}{\simeq }M^{\circ }$. Let $L\subset
M $ be an $R$-cofinite $A$-submodule with $\{f_{1},...,f_{k}\}$ a generating
system of $\mathrm{An}(L)\simeq (M/L)^{\ast }.$ Then there exists by
assumption $\{n_{1},...,n_{k}\}\subset N,$ such that $\chi
_{Q}(n_{i})=f_{i}. $ By \cite[28.1]{Wis88}$\ $we have then 
\begin{equation*}
(\sum\limits_{i=1}^{k}Rn_{i})^{\bot }=\bigcap\limits_{i=1}^{k}\mathrm{Ke}%
(f_{i})=\mathrm{Ke}(\sum\limits_{i=1}^{k}Rf_{i})=\mathrm{KeAn}(L)=L,
\end{equation*}
i.e. $L$ is open w.r.t. $M[\frak{T}_{ls}(N)].$ Consequently $\mathrm{Cf}%
(M)\leq M[\frak{T}_{ls}(N)].$ By Lemma \ref{leq} (2) $M[\frak{T}%
_{ls}(N)]\leq \mathrm{Cf}(M)$ and so $M[\frak{T}_{ls}(N)]=\mathrm{Cf}(M),$
i.e. $Q$ is coreflexive.

\item  The result follows from (1) and (a).$\blacksquare $
\end{enumerate}
\end{enumerate}
\end{Beweis}

\begin{corollary}
\label{N-koref}Let $R$ be noetherian and $C$ be a locally projective $R$%
-coalgebra.

\begin{enumerate}
\item  If $N$ is coreflexive, then $N^{\ast \Box }=N^{\ast \circ }.$

\item  Let $N^{\ast }$ be cofinitely $R$-cogenerated.

\begin{enumerate}
\item  If $N\simeq N^{\ast \circ },$ then $N$ is coreflexive.

\item  Let $N$ be weakly coreflexive. Then $N$ is coreflexive, iff $N\simeq
N^{\ast \circ }.$
\end{enumerate}
\end{enumerate}
\end{corollary}

\begin{theorem}
\label{cor}Let $R$ be noetherian and $P=(A,C)\in \mathcal{P}_{m}^{\alpha }$.

\begin{enumerate}
\item  If $P$ is coreflexive, then $C=A^{\circ }.$

\item  Assume $R$ to be artinian. Then $P$ is coreflexive, iff all $R$%
-cofinite $A$-ideals are closed w.r.t. $A[\frak{T}_{ls}(C)]=\mathcal{T}%
_{-C}(A).$

\item  If $A$ is cofinitely $R$-cogenerated, then the following statements
are equivalent:

(i) $P$ is coreflexive;

(ii) $C=A^{\circ }$.

(iii) every locally finite left $A$-module is $C$-rational, i.e. $\mathrm{Loc%
}(_{A}\mathcal{M})=\sigma \lbrack _{A}C]${\normalsize .}
\end{enumerate}
\end{theorem}

\begin{Beweis}
\begin{enumerate}
\item  By Corollary \ref{C=A^r} (1) $C=A^{r}$ and so the result follows from
Theorem \ref{Q-koref} (1).

\item  Let $R$ be artinian. By \cite[Lemma 1.7 (4)]{Abu} every $R$-cofinite
closed $A$-ideal is open and the result follows.

\item  (i) $\Leftrightarrow $ (ii)\ follows from Theorem \ref{Q-koref} (3).

(ii)\ $\Rightarrow $ (iii) By assumption and Proposition \ref{A^0-cor} (2) $%
\mathrm{Loc}(_{A}\mathcal{M})=\sigma \lbrack _{A}A^{\circ }]=\sigma \lbrack
_{A}C].$

(iii) $\Rightarrow $ (ii) Assume all locally finite left $A$-modules to be $%
C $-rational. Then in particular $_{A}A^{\circ }$ is $C$-rational and it
follows from Corollary \ref{C=A^r} (2) that $C=A^{\circ }.\blacksquare $
\end{enumerate}
\end{Beweis}

\begin{corollary}
\label{C-coref}Let $R$ be noetherian and $C$ be a locally projective $R$%
-coalgebra.

\begin{enumerate}
\item  If $C$ is coreflexive, then the canonical $R$-linear mapping $\phi
_{C}:C\rightarrow C^{\ast \ast }$ induces an isomorphism $C\overset{\Phi _{C}%
}{\simeq }C^{\ast \circ }.$

\item  Let $R$ be artinian. Then $C$ is coreflexive, iff all $R$-cofinite $%
C^{\ast }$-ideals are closed w.r.t. the finite topology.

\item  If $C^{\ast }$ is cofinitely $R$-cogenerated, then the following
statements are equivalent:

(i) $C$ is coreflexive;

(ii) $C\simeq C^{\ast \circ };$

(iii) every locally finite left $C^{\ast }$-module is $C$-rational.
\end{enumerate}
\end{corollary}

\qquad As a consequence of Lemma \ref{grun} and Theorem \ref{cor} (3) get we

\begin{proposition}
\label{M0-char}Let $R$ be noetherian. If $A$ is a cofinitely $R$-cogenerated 
$\alpha $-algebra and $M$ is a right $A$-module with structure map $\phi
_{M}:M\otimes _{R}A\rightarrow M,$ then for every $f\in M^{\ast }$ the
following statements are equivalent:

\begin{enumerate}
\item  $f\in M^{\circ }.$

\item  $\phi _{M}^{\ast }(f)\in M^{\circ }\otimes _{R}A^{\circ }.$

\item  $\phi _{M}^{\ast }(f)\in M^{\circ }\otimes _{R}A^{\ast }.$

\item  $Af$ is finitely generated in $\mathcal{M}_{R}.$

\item  $f(MI)=0$ for an $R$-cofinite \emph{(}right\emph{)} $A$-ideal.

\item  $f(L)=0$ for an $R$-cofinite right $A$-submodule $L\subset M.$
\end{enumerate}
\end{proposition}

Analog to \cite{Taf72} we get

\begin{corollary}
\label{ref}Let $R$ be a QF ring.

\begin{enumerate}
\item  A projective $R$-coalgebra $C$ is coreflexive, iff $C^{\ast }$ is a
reflexive $R$-algebra.

\item  Let $A$ be an $\alpha $-algebra. If $A$ is weakly reflexive, then $%
A^{\circ }$ is a coreflexive $R$-coalgebra.
\end{enumerate}
\end{corollary}

\begin{example}
\label{abz}(\cite[Example 5]{Lin77}) Let $R$ be a field and consider the
Hopf $R$-algebra $(H,\mu ,\eta ,\Delta ,\varepsilon ,S)$ with countable
basis $\{h_{0},h_{1},h_{2},...\}$ and 
\begin{equation*}
\begin{tabular}{lllllllll}
$\mu (h_{n}\otimes h_{k})$ & $:=$ & $\binom{n+k}{n}h_{n+k},$ & $\Delta
(h_{n})$ & $:=$ & $\sum\limits_{i+j=n}h_{i}\otimes h_{j},$ & $S(h_{n})$ & $%
:= $ & $(-1)^{n}h_{n}.$ \\ 
$\eta (1_{R})$ & $:=$ & $h_{0},$ & $\varepsilon (h_{n})$ & $:=$ & $\delta
_{0,n}.$ &  &  & 
\end{tabular}
\end{equation*}

\begin{enumerate}
\item  $H^{\ast }\simeq R[[x]]$ is a principal ideal domain and 
\begin{equation*}
\mathcal{M}{\normalsize ^{H}\simeq \mathrm{Rat}^{H}(_{H^{\ast }}}\mathcal{M}%
{\normalsize )=\{M\in }\text{ }{\normalsize _{H^{\ast }}}\mathcal{M}%
{\normalsize |}\text{ }M\text{ is a torsion module}\}{\normalsize .}
\end{equation*}

So $\mathrm{Rat}^{H}(-)$ is a radical and $\mathrm{Rat}^{H}(_{H^{\ast }}%
\mathcal{M})$ is closed under extensions.

\item  $H^{\Box }:=\mathrm{Rat}^{H}(_{H}H^{\ast })=0.$

\item  There exists no finite dimensional nonzero projective right $H$%
-comodules.

\item  $H\simeq H^{\ast \circ },$ i.e. $H$ is a coreflexive $R$-coalgebra.
\end{enumerate}
\end{example}

\textbf{Acknowledgment.} This paper (up to a few changes) includes parts of
my doctoral thesis at the Heinrich-Heine Universit\"{a}t, D\"{u}sseldorf
(Germany). I am so grateful to my advisor Prof. Robert Wisbauer for his
wonderful supervision and the continuous encouragement and support.

\end{document}